\documentclass[a4paper,12pt]{amsart}

\usepackage{amsmath,times,epsfig,amssymb,amsbsy,amscd,amsfonts,amstext,color,bm}

\usepackage[all]{xy}
\usepackage{bm}
\usepackage{comment}
\usepackage{enumerate}
\usepackage{circledsteps}

\usepackage{tikz,epic,eso-pic}

\tikzset{v/.style={circle, draw, inner sep=2pt, minimum size=6pt, fill=white}}

\theoremstyle{plain}
\newtheorem{theorem}{Theorem}[section]
\newtheorem*{theorem*}{Theorem}
\newtheorem*{theoremA*}{Theorem A}
\newtheorem*{theoremB*}{Theorem B}
\newtheorem{corollary}[theorem]{Corollary}
\newtheorem{lemma}[theorem]{Lemma}
\newtheorem{conjecture}[theorem]{Conjecture}
\newtheorem{problem}[theorem]{Problem}

\theoremstyle{definition}

\newtheorem{remark}[theorem]{Remark}
\newtheorem{example}[theorem]{Example}
\newtheorem{proposition}[theorem]{Proposition}

\newcommand{\corank}{\operatorname{corank}}

\newcommand{\bC}{\mathbb{C}}

\newcommand{\bZ}{\mathbb{Z}}

\newcommand{\cA}{\mathcal{A}}

\newcommand{\cC}{\mathcal{C}}

\def\qed{\hfill $\Box$}

\title{First homology groups of the Milnor fiber boundary for generic hyperplane  arrangements in $\bC^{3}$}
\author{Sakumi Sugawara}
\date{\today}
\address{Department of Mathematics, Graduate School of Science, Hokkaido University, North 10, West 8, Kita-ku, Sapporo 060-0810, JAPAN. }
\email{sugawaras@math.sci.hokudai.ac.jp}

\begin{document}
\maketitle

\begin{abstract}
We study the Milnor fiber boundary for hyperplane arrangements in $\bC^3$.
This is one of the examples of non-isolated surface singularities, which are studied by N\'emethi--Szil\'ard.
In this paper, we compute the first homology group of the Milnor fiber boundary for a generic arrangement, which gives an affirmative answer to the conjecture of Suciu.
Also, we give an example of an arrangement with $n$ hyperplanes, whose torsion part in the Milnor fiber boundary homology contains a direct summand other than $\bZ_{n}$, for certain value of $n$. 
\end{abstract}


\section{Introduction}
Let $(V,\mathbf{0})$ be a germ of an analytic hypersurface in a complex Euclidean space with a singularity at the origin.
The link of an isolated singularity is defined as its transverse intersection with a small sphere centered at the origin.
The singularity link is a very important object in singularity theory and topology.
Starting with Milnor's fibration theorem \cite{mil-sing}, it has been studied by a lot of authors for a long time \cite{sea}.
It is known that the singularity link of a complex surface with an isolated singularity is a $3$-manifold obtained by plumbing along the dual graph of the minimal resolution \cite{mum, neu}.
If the singularity is not isolated, the topological description of the singularity may be complicated in general.
For a non-isolated surface singularity, N\'emethi--Szil\'ard and Curmi showed that the boundary of the Milnor fiber, which is the counterpart of the link in the isolated singularity, is a plumbed $3$-manifold.
They also described the detailed construction of the plumbing graph \cite{nem-szi, cur}.
However, there are still a lot of open problems with the topology of the Milnor fiber boundary of a non-isolated surface singularity (see Section 24.4 in \cite{nem-szi}).

Among non-isolated surface singularities, we will study a hyperplane arrangement in $\bC^3$ in this paper.
A hyperplane arrangement is a finite set of linear hyperplanes in a vector space.
The union of hyperplanes defines a variety with a non-isolated singularity at the origin.
Hyperplane arrangements are studied from various aspects, for example, combinatorics, topology, and algebraic geometry \cite{orl-ter}.
From a topological viewpoint, one of the important problems is whether several topological invariants are combinatorially determined (i.e., determined by the information of the intersection poset).
For example, the cohomology ring of the complement is combinatorially determined \cite{orl-sol}, while the fundamental group is not \cite{ryb}.
There are many open problems about whether topological invariants of covering spaces and Milnor fibers are combinatorially determined, and even for the first Betti number (see \cite{pap-suc} for recent progress on this topic).

In this paper, we study the Milnor fiber boundary of hyperplane arrangements in $\bC^3$.
The Milnor fiber boundary of a hyperplane arrangement has been studied in \cite{nem-szi, suc-mil, suc-mil2} so far. 
N\'emethi--Szilard described the Milnor fiber boundary as a plumbed $3$-manifold whose plumbing graph is defined by its combinatorics (Subsection 8.2 in \cite{nem-szi} and Section \ref{sec:algo} in this paper). 
Thus, the homeomorphism type of the Milnor fiber boundary is combinatorially determined.
Moreover, they described a formula for the first Betti number of the Milnor fiber boundary (Theorem 7.8 in \cite{suc-mil}, also Theorem 19.10.2 in \cite{nem-szi}).
\begin{theorem} 
Let $\cA$ be a hyperplane arrangement in $\bC^3$ with $n$ hyperplanes and denote the Milnor fiber boundary by $\partial \overline{F}$. The first Betti number of $\partial \overline{F} $ is expressed as:
\[
b_{1} (\partial \overline{F}) = \sum_{X \in L_{2} (\cA)} (1 + (|\cA_{X}|-2) \gcd(|\cA_{X}|, n) ).
\]
(The definitions of notations will be given in Section \ref{sec:arr}.)
\end{theorem}

However, the combinatorial formula for the integral first homology group is still unknown. 
The following can be read from Problem 24.4.19 in \cite{nem-szi}:
\begin{problem}
Find a nice formula for the torsion of $H_{1} (\partial \overline{F}, \bZ)$.
\end{problem}
In general, the Milnor fiber boundary may have a torsion in the first homology group. 
In fact, it is known that $H_{1} (\partial \overline{F} , \bZ) = \bZ^{6} \oplus \bZ_4$ for a generic arrangement with four hyperplanes.
In this paper, we attack the above problem somewhat.
We compute the first integral homology group of the Milnor fiber boundary for a generic arrangement. 
The following is the main result.
This gives the affirmative answer to the conjecture of Suciu (page 468 in \cite{suc-mil}).
\begin{theorem}
(=Theorem \ref{thm:main}.)
For a generic hyperplane arrangement with $n$ hyperplanes, the first homology group of the Milnor fiber boundary $\partial \overline{F}$ is expressed as:
\[
H_{1} (\partial \overline{F}, \bZ)  = \bZ^{n(n-1)/2} \oplus \bZ_{n}^{(n-2)(n-3)/2}.
\]
\end{theorem}
It is known that the first homology group of a plumbed manifold is computed by the weighted incidence matrix of the plumbing graph (see Theorem \ref{thm:homology}).
Following the book \cite{nem-szi}, we can construct the plumbing graph of the Milnor fiber boundary of arrangements. 
By computing the incidence matrix and its Smith normal form, we obtain the first homology of the Milnor fiber boundary.

Using their technique, we can compute the first homology groups of the Milnor fiber boundary for any arrangements.
The author's examination so far deduces the following conjecture.

\begin{conjecture}\label{conj:tor}
Let $\cA$ be a central hyperplane arrangement in $\bC^{3}$ with $n$ hyperplanes and $\partial \overline{F}$ be its Milnor fiber boundary.
\begin{enumerate}[(1)]
\item Suppose that $(|\cA_{X}|-2)(\gcd (|\cA_{X}|, n)-1)=0$
for each $X \in L_{2} (\cA)$. Then the torsion of $H_{1} (\partial \overline{F}; \bZ)$ is the direct sum of $\bZ_{n}$ with the number equals to the Euler characteristic $\chi(\bC P^2 \setminus \bigcup_{H \in \cA} \overline{H})$ of the projectivized complement.
In particular, the torsion part is the direct sum of $\bZ_{n}$ for any arrangement if $n$ is prime.
\item If $\bZ_{m}$ is contained in the torsion of $H_{1} (\partial \overline{F}; \bZ)$, then $m|n$. 
\item $H_{1} (\partial \overline{F}; \bZ)$ is torsion free if and only if $\cA$ is either pencil or near pencil.
\end{enumerate}
\end{conjecture} 
Note that if the assumption of (1) fails, then there is an arrangement with $n$ hyperplanes, whose torsion part in $H_{1}(\partial F; \bZ)$ contains a direct summand other than $\bZ_{n}$, see Example \ref{exam:a3}.
This gives a counterexample for an open question written on page 210 in \cite{suc-mil2}, which asks the torsion part in $H_{1} (\partial F; \bZ)$ for any arrangement with $n$ hyperplanes.

This paper is organized as follows. 
First, we review the basic notations on hyperplane arrangements in Section \ref{sec:arr}, and plumbed manifolds in Section \ref{sec:graph}.
We give the proof of the main result in Section \ref{sec:main}.
Finally, we give the algorithm to obtain the plumbing graph of the Milnor fiber boundary specialized in the case of arrangements in Section \ref{sec:algo}.

\vspace{3mm}
\textbf{Acknowledgements.}
The author would like to thank Professor Naohiko Kasuya for the helpful discussions and comments on this paper. 
The author also thanks Masato Tanabe for the discussion on Lemma \ref{lem:keylem},
and Professor Alexander Suciu for the encouragement on this work.
The author also thanks Yongqiang Liu for the discussion on Conjecture \ref{conj:tor}.
This work is supported by JSPS KAKENHI Grant Number 22KJ0114.

\section{Hyperplane arrangements}\label{sec:arr}

A finite collection $\cA = \{H_1 , \ldots, H_{n} \}$ of linear hyperplanes in a complex vector space $\bC^{\ell + 1}$ is called a \textit{hyperplane arrangement}.
We set $M(\cA) = \bC^{\ell + 1} \setminus \bigcup_{i=1}^{n} H_{i}$.
The set of non-empty intersections is called the \textit{intersection poset} and is denoted by $L(\cA)$. 
We denote the set of intersections of codimension $k$ by $L_{k}(\cA)$. 
For example, $L_{0}(\cA) = \{\bC^{\ell+1}\}$ and $L_{1}(\cA) = \cA$.
For an intersection $X \in L(\cA)$, the localization to $X$ is denoted by $\cA_{X}$, that is, $\cA_{X} = \{H \in \cA \mid H \supset X\}$.
An arrangement $\cA$ is called \textit{generic} if $|\cA_{X}| = 2$ for all $X \in L_{2} (\cA)$.

Let $\alpha_{i}$ be a defining linear form of $H_{i}$.
The product $Q = \prod_{i=1}^{n} \alpha_{i}$ defines a polynomial function $Q: \bC^{\ell + 1} \rightarrow \bC$.
The restriction to the non-zero locus $Q: M(\cA) \rightarrow \bC^{*}$ is known to be a fiber bundle \cite{mil-sing}.
Note that the singularity of $Q^{-1} (0) = \bigcup_{i=1}^{n} H_{i}$ is non-isolated if $\ell \geq 2$.
The regular fiber $F = Q^{-1} (1)$ is called the \textit{Milnor fiber} of $\cA$.
For a sufficiently large ball $B^{2\ell+2} \subset \bC^{\ell+1}$, the intersection $\overline{F}= F \cap B^{2\ell+2}$ is called the closed Milnor fiber. 
The boundary $\partial \overline{F} = F \cap \partial B^{2\ell+2}$ of the closed Milnor fiber is called the \textit{Milnor fiber boundary}. 
The Milnor fiber boundary is a smooth closed $(2\ell-1)$-dimensional manifold.
We will deal with $\partial \overline{F}$ as the main topic in this paper.

Each linear form $\alpha_{i}$ defines a projective hyperplane $\overline{H_{i}} \subset \bC P^{\ell}$.
Let $\overline{\cA} = \{\overline{H_{1}} , \ldots, \overline{H_{n}}\}$ be the projectivized arrangement in $\bC P^{\ell}$. 
Let us consider the case where $\ell=2$, that is, $\overline{\cA}$ is a line arrangement in $\bC P^2$.
Let $P(\overline{\cA}) = \{p_{1} ,\ldots, p_{t}\}$ be the set of intersection points of $\overline{\cA}$.
We define a bipartite graph $\Gamma (\overline{\cA})$ as follows:
\begin{itemize}
\item For $\overline{H_{i}} \in \overline{\cA}$, we give a corresponding vertex $v_{i}$ (line vertex).
\item For $p_{j} \in P(\overline{\cA})$, we give a corresponding vertex $w_{j}$ (point vertex).
\item Vertices $v_{i}$ and $w_{j}$ are connected by an edge if and only if  $p_{j} \in \overline{H_{i}}$.
\end{itemize}
The graph $\Gamma (\overline{\cA})$ is called the \textit{incidence graph} of $\overline{\cA}$.
Note that $L(\cA)$ and $\Gamma(\overline{\cA})$ have equivalent combinatorial information.

\section{Plumbed manifolds}\label{sec:graph}
In this section, we review plumbed $3$-manifolds (see \cite{neu} and Chapter 4 in \cite{nem-szi} for details).

\subsection{Plumbing graphs and plumbing calculus}

Let $\Gamma$ be a finite graph. We denote the vertex set by $V(\Gamma)$ and the edge set by $E(\Gamma)$. 
An \textit{oriented closed plumbing graph} is the graph with the following decoration:
\begin{itemize}
\item Each vertex $i$ has two kinds of weight; $g_{i} \in \bZ_{\geq 0}$ (genus) and $e_{i} \in \bZ$ (Euler number).
\item Each edge equips a signature $+$ or $-$.
\end{itemize}
The weight of the Euler number is sometimes just called weight.
The genus of a vertex is written as $[g_{i}]$ and omitted if it is zero.
The signature of an edge $e$ is written as $\varepsilon_{e}$.

From an oriented closed plumbing graph (shortly, plumbing graph), we can construct a closed $3$-manifold by the following procedure:
Corresponding to a vertex $i$, we give an oriented $S^1$-bundle $\pi_{i} : E_{i} \rightarrow S_{i}$ over a genus $g_{i}$ orientable surface $S_{i}$ with the Euler number $e_{i}$.
Let $p \in S_{i}$ and $D_{p}$ be a small disk centered at $p$.
Then, we can locally trivialize the fibration as $\pi_{i}^{-1} (D_{p}) \cong D_{p} \times S^1$.
Suppose that a vertex $j$ is connected with $i$ by an edge of signature $\varepsilon \in \{\pm \}$. 
Take a point $q \in S_{j}$ and trivialize around its neighborhood $\pi_{j}^{-1} (D_{q}) \cong D_{q} \times S^1$, similarly. 
Then, glue $\pi_{i}^{-1} (S_{i}\setminus D_{p})$ and $\pi_{j}^{-1} (S_{j} \setminus D_{q})$ along a homeomophism $\varphi : \partial D_{p} \times S^1 \rightarrow \partial D_{q} \times S^1$ which is presented by a matrix 
$\varepsilon \bigl(
\begin{smallmatrix}
   0 & 1 \\
   1 & 0
\end{smallmatrix}
\bigl)$.
After repeating this operation for every vertex and edge, we obtain a closed $3$-manifold $M = M(\Gamma)$ called a \textit{plumbed $3$-manifold}, or a \textit{graph $3$-manifold}.

There are some moves for a plumbing graph that do not change the topology of the resulting plumbed manifold, called \textit{plumbing calculus}.
We will list only plumbing calculus to be used in this paper (see \cite{neu} for the complete list).

\begin{proposition}
Applying the following operations to a plumbing graph does not change the diffeomorphism type of $M(\Gamma)$:
\begin{itemize}
\item [(i)] Reversing the signature of all non-loop edges adjacent to a fixed vertex.
\item [(ii)] \textbf{Blowing down.} Here, $\varepsilon = \pm 1$ and the signatures satisfy $\varepsilon_{0} = - \varepsilon \varepsilon_{1} \varepsilon_{2}$.
	\begin{enumerate}[(a)]
	\item 
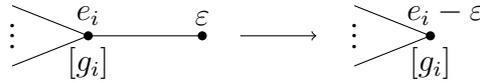
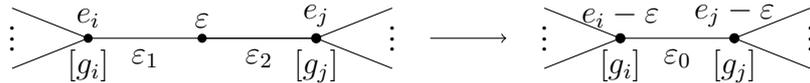
	\[
	\begin{tikzpicture}
	\fill (0,0) circle (0.06);
	\draw (0,0) node[above]{$e_{i}$};
	\draw (0,0) node[below]{$[g_{i}]$};
	\draw (0,0) --++(-1,0.4);
	\draw (0,0) --++(-1,-0.4);
	\draw (-1,0.1)node{$\vdots$};
	\draw (0,0) --++(1.5,0);
	\draw (1.5,0) node[above]{$\varepsilon$};
	\fill (1.5,0) circle(0.06);
	
	\draw[->] (2,0) --++(1,0);
	\fill (4.5,0) circle (0.06);
	\draw (4.7,0) node[above]{$e_{i}-\varepsilon$};
	\draw (4.5,0) node[below]{$[g_{i}]$};
	\draw (4.5,0) --++(-1,0.4);
	\draw (4.5,0) --++(-1,-0.4);
	\draw (3.5,0.1)node{$\vdots$};
	\end{tikzpicture}
	\]
	\item 
	\[
	\begin{tikzpicture}
	\fill (0,0) circle (0.06);
	\draw (0,0) node[above]{$e_{i}$};
	\draw (0,0) node[below]{$[g_{i}]$};
	\draw (0,0) --++(-1,0.4);
	\draw (0,0) --++(-1,-0.4);
	\draw (-1,0.1)node{$\vdots$};
	\draw (0,0) --++(1.5,0);
	\draw (1.5,0) node[above]{$\varepsilon$};
	\fill (1.5,0) circle(0.06);
	\draw (0.75,0) node[below]{$\varepsilon_{1}$};
	\draw (2.25,0) node[below]{$\varepsilon_{2}$};	
	
	\draw (1.5,0)--++(1.5,0);
	\filldraw (3,0) circle (0.06) node[above]{$e_{j}$};
	\draw (3,0)node[below]{$[g_{j}]$};
	\draw (3,0) --++(1,0.4);
	\draw (3,0) --++(1,-0.4);
	\draw (4,0.1)node{$\vdots$};
	
	\draw[->] (4.5,0) --++(1,0);
	
	\filldraw (7,0) circle (0.06) node[above]{$e_{i}-\varepsilon$};
	\draw (7,0) node[below]{$[g_{i}]$};
	\draw (7,0) --++(-1,0.4);
	\draw (7,0) --++(-1,-0.4);
	\draw (6,0.1)node{$\vdots$};
	\draw (7,0) --++(1.5,0);
	
	\draw (7.75,0) node[below]{$\varepsilon_{0}$};	
	
	\draw (1.5,0)--++(1.5,0);
	\filldraw (8.5,0) circle (0.06) node[above]{$e_{j}-\varepsilon$};
	\draw (8.5,0)node[below]{$[g_{j}]$};
	\draw (8.5,0) --++(1,0.4);
	\draw (8.5,0) --++(1,-0.4);
	\draw (9.5,0.1)node{$\vdots$};
	
	\end{tikzpicture}
	\]

	\end{enumerate}
\item [(iii)] \textbf{$0$-chain absorption.} Each signature satisfies $\varepsilon'_{i} = -\varepsilon \overline{\varepsilon} \varepsilon_{i}$ if the edge is not a loop, and $\varepsilon'_{i} = \varepsilon_{i}$ if it is a loop.
	
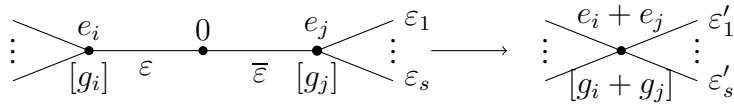
\[
	\begin{tikzpicture}
	\fill (0,0) circle (0.06);
	\draw (0,0) node[above]{$e_{i}$};
	\draw (0,0) node[below]{$[g_{i}]$};
	\draw (0,0) --++(-1,0.4);
	\draw (0,0) --++(-1,-0.4);
	\draw (-1,0.1)node{$\vdots$};
	\draw (0,0) --++(1.5,0);
	\draw (1.5,0) node[above]{$0$};
	\fill (1.5,0) circle(0.06);
	\draw (0.75,0) node[below]{$\varepsilon$};
	\draw (2.25,0) node[below]{$\overline{\varepsilon}$};	
	
	\draw (1.5,0)--++(1.5,0);
	\filldraw (3,0) circle (0.06) node[above]{$e_{j}$};
	\draw (3,0)node[below]{$[g_{j}]$};
	\draw (3,0) --++(1,0.4) node[right]{$\varepsilon_{1}$};
	\draw (3,0) --++(1,-0.4) node[right]{$\varepsilon_{s}$};
	\draw (4,0.1)node{$\vdots$};	
	
	\draw[->] (4.5,0) --++(1,0);
	\draw (7,0)--++(-1,0.4);
	\draw (7,0)--++(-1,-0.4);
	\draw (6,0.1)node{$\vdots$};
	\draw (7,0)--++(1,0.4) node[right]{$\varepsilon'_{1}$};
	\draw (7,0)--++(1,-0.4) node[right]{$\varepsilon'_{s}$};
	\draw (8,0.1)node{$\vdots$};
	\fill (7,0) circle (0.06) ;
	\draw (7,0.1)node[above]{$e_{i} + e_{j}$};
	\draw (7,-0.1) node[below]{$[g_{i} + g_{j}]$};
	\end{tikzpicture}
	\]
\item [(iv)] \textbf{Oriented handle absorption.} 
	
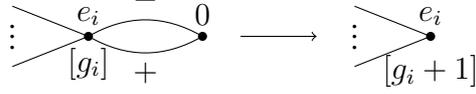
\[
	\begin{tikzpicture}
	\fill (0,0) circle (0.06);
	\draw (0,0) node[above]{$e_{i}$};
	\draw (0,0) node[below]{$[g_{i}]$};
	\draw (0,0) --++(-1,0.4);
	\draw (0,0) --++(-1,-0.4);
	\draw (-1,0.1)node{$\vdots$};
	\draw (0,0) to[out=30,in=150] ++(1.5,0);
	\draw (0.75,0.2) node[above]{$-$};
	\draw (0,0) to[out=-30,in=-150] ++(1.5,0);
	\draw (0.75,-0.2) node[below]{$+$};
	\fill (1.5,0) circle (0.06);
	\draw (1.5,0) node[above]{$0$};
	
	\draw[->] (2,0) --++ (1,0);
	\fill (4.5,0) circle (0.06);
	\draw (4.5,0) node[above]{$e_{i}$};
	\draw (4.5,-0.1) node[below]{$[g_{i}+1]$};
	\draw (4.5,0) --++(-1,0.4);
	\draw (4.5,0) --++(-1,-0.4);
	\draw (3.5,0.1)node{$\vdots$};
	\end{tikzpicture}
	\]
\item [(v)] \textbf{Splitting.} Each $\Gamma_{j}$ is a connected graph and vertex $i$ is connected to $\Gamma_{j}$ by $k_{j}$ edges. The graph decomposes into the disjoint union of $\Gamma_{1} , \ldots , \Gamma_{t}$ and $(2g_{i} + \sum_{j}(k_{j}-1))$ copies of a vertex with weight $0$.
	
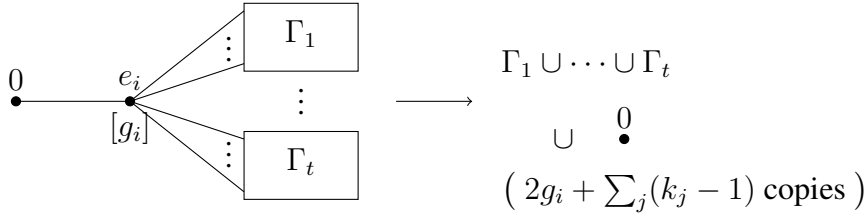
\[
	\begin{tikzpicture}
	\filldraw (0,0) circle (0.06) node[above]{$0$};
	\draw (0,0) --++(1.5,0);
	\filldraw (1.5,0) circle (0.06) node[above]{$e_{i}$};
	\draw (1.5,0) node[below]{$[g_{i}]$};
	
	\draw (1.5,0) --++(1.5,1.2);
	\draw (1.5,0) --++(1.5,0.5);
	\draw (2.8,0.8) node{$\vdots$};
	\draw (1.5,0) --++(1.5,-0.5);
	\draw (1.5,0) --++(1.5,-1.2);
	\draw (2.8,-0.6) node{$\vdots$};
	
	\draw (3,1.3) --++(1.5,0)--++(0,-0.9)--++(-1.5,0)--cycle;
	\draw (3.75,0.9) node{$\Gamma_{1}$};
	\draw (3.75,0.1) node{$\vdots$};
	\draw (3,-0.4) --++(1.5,0)--++(0,-0.9)--++(-1.5,0)--cycle;
	\draw (3.75,-0.8) node{$\Gamma_{t}$};
	
	\draw[->] (5,0) --++(1,0);
	\draw (7.5,0.5)node{$\Gamma_{1} \cup \cdots \cup \Gamma_{t}$};
	\draw (7.2,-0.5) node{$\cup$};
	\filldraw (8,-0.5) circle (0.06) node[above] {$0$};
	\draw (8.8,-1.2) node{$\big ($ $2g_{i} + \sum_{j}(k_{j}-1)$ copies $\big )$};
	\end{tikzpicture}
	\]
\end{itemize}
\end{proposition}

The following calculus will be useful later.
We call this \textit{$\pm 2$-alteration}.
It is obtained by performing a $(-1)$-blow-up, and then a $1$-blow-down.
\begin{itemize}
\item[(b')] \textbf{$\pm 2$-alteration}.
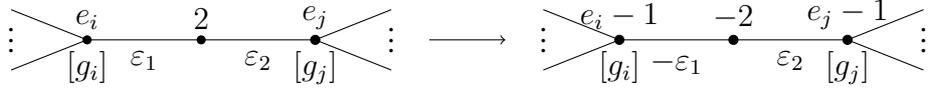
\[
	\begin{tikzpicture}
	\fill (0,0) circle (0.06);
	\draw (0,0) node[above]{$e_{i}$};
	\draw (0,0) node[below]{$[g_{i}]$};
	\draw (0,0) --++(-1,0.4);
	\draw (0,0) --++(-1,-0.4);
	\draw (-1,0.1)node{$\vdots$};
	\draw (0,0) --++(1.5,0);
	\draw (1.5,0) node[above]{$2$};
	\fill (1.5,0) circle(0.06);
	\draw (0.75,0) node[below]{$\varepsilon_{1}$};
	\draw (2.25,0) node[below]{$\varepsilon_{2}$};	
	
	\draw (1.5,0)--++(1.5,0);
	\filldraw (3,0) circle (0.06) node[above]{$e_{j}$};
	\draw (3,0)node[below]{$[g_{j}]$};
	\draw (3,0) --++(1,0.4);
	\draw (3,0) --++(1,-0.4);
	\draw (4,0.1)node{$\vdots$};
	
	\draw[->] (4.5,0) --++(1,0);
	
	\filldraw (7,0) circle (0.06) node[above]{$e_{i}-1$};
	\draw (7,0) node[below]{$[g_{i}]$};
	\draw (7,0) --++(-1,0.4);
	\draw (7,0) --++(-1,-0.4);
	\draw (6,0.1)node{$\vdots$};
	\draw (7,0) --++(1.5,0);
	
	\filldraw (8.5,0) circle (0.06) node[above]{$-2$};
	\draw (7.75,0) node[below]{$-\varepsilon_{1}$};
	\draw (8.5,0) --++(1.5,0);	
	\draw (9.25,0) node[below]{$\varepsilon_{2}$};	
	\filldraw (10,0) circle (0.06) node[above]{$e_{j}-1$};
	\draw (10,0)node[below]{$[g_{j}]$};
	\draw (10,0) --++(1,0.4);
	\draw (10,0) --++(1,-0.4);
	\draw (11,0.1)node{$\vdots$};
	
	\end{tikzpicture}
	\]

\end{itemize}

\subsection{Multiplicity systems}
We can give extra weight called multiplicity to each vertex (see Section 4.1.4 in \cite{nem-szi} for details).
From a plumbing graph $\Gamma$, we can construct the plumbed $4$-manifold in a similar way as the plumbed $3$-manifold.
Corresponding to each vertex, we give the associated $D^2$-bundle. 
By gluing them together in a similar way, we obtain a $4$-manifold $P=P(\Gamma)$.
Clearly, $\partial P = M$.
For each vertex $i$, the zero section of the $D^2$-bundle represents a $2$-cycle $C_{i} \in H_{2} (P;\bZ)$.
Under this setting, the multiplicity system $\{m_{i}\}_{i \in V}$ is is a family of intergers $m_{i}$ which satisfy the relation $\sum_{i \in V} m_{i} C_{i} = 0 \in H_{2} (P, \partial P; \bZ)$.
The multiplicity is written as $(m_{i})$ in the graph.

It is known that the multiplicity satisfies the following local formula for each vertex:
\begin{equation}\label{eq:mult}
e_{i} m_{i} + \sum_{e \in E_{i}} \varepsilon_{e}  m_{v(e)} =0,
\end{equation}
where $E_{i}$ is the set of edges adjacent to $i$ and the edge $e$ connects vertices $i$ and $v(e)$.
Thus, once we know the multiplicity system of a plumbing graph, we can compute the Euler number of each vertex from this formula.
Note that when computing the multiplicity of each vertex, the plumbing graph can have arrowheads, and an arrowhead has a multiplicity. 
However, these arrowheads are not used for the definition for closed plumbed manifolds.
In Section \ref{sec:algo}, we will determine the plumbing graph of a Milnor fiber boundary. At that time, the multiple system of the graph is given first. 
Thus, we use the above formula and then compute Euler numbers.

\subsection{The first homology group}
Next, we will explain how to compute the first homology group for a plumbed $3$-manifold.
For simplicity, we assume that $\Gamma$ is a simple graph, that is, $\Gamma$ does not have loops and multiple edges.
Let $g(\Gamma) = \sum_{i \in V} g_{i}$ and put $W$ by the set of non-arrowhead vertices.
Let denote $A=(a_{ij})_{i,j \in W}$ the weighted incidence matrix of $\Gamma$, that is, each entry is defined as
\begin{eqnarray*}
a_{ij} = 
\left \{
\begin{array}{ll}
e_{i} & (i=j), \\
\varepsilon_{(i,j)} & (i \neq j).
\end{array}
\right .
\end{eqnarray*}
The first homology group of the plumbed manifold is computed as follows (see Chapter 15 in \cite{nem-szi}, also Lemma 3.4 in \cite{acm}). 
\begin{theorem}\label{thm:homology}
The first Betti number of $M$ is computed as 
\[
b_{1} (M) = \corank(A) + 2 g(\Gamma) + b_{1} (\Gamma),
\]
and the torsion of $H_{1} (M;\bZ)$ coincides with the torsion of the cokernel of a $\bZ$-module homomorphism represented by the matrix $A$.
\end{theorem}
Note that the torsion part is determined by the Smith normal form of $A$.

\section{Computation of the homology group}\label{sec:main}
Throughout this section, We assume that $\cA$ is a hyperplane arrangement with $n$ generic hyperplanes defined in $\bC^{3}$ and denote its Milnor fiber boundary by $\partial \overline{F}$. 

As we will explain in Section \ref{sec:algo} (see also \cite{nem-szi}), $\partial \overline{F}$ is a plumbed $3$-manifold.
The plumbing graph of $\partial \overline{F}$ is obtained as follows: 
First, we consider the complete graph of $n$ vertices and give a weight $-1$ to each vertex.
Next, insert a vertex with weight $-n$ to each edge.
Each edge of the complete graph is divided into exactly two edges, one is decorated by $+$ and the other is decorated by $-$. 
The genus is zero for each vertex (see Example \ref{exam:main} for the concrete construction of the plumbing graph).

It is convenient to give a total order to the vertex set. 
First, order the vertices of the original complete graph as $v_{1}, \ldots, v_{n}$. 
Let $w_{ij}$ be the vertex inserted into the edge connecting $v_{i}$ and $v_{j}$ ($i<j$). 
Then, order vertices as follows:
\begin{itemize}
\item $v_{i} < v_{j}$ \; for each $i<j$.
\item $v_{i} < w_{jk}$ \; for each $i,j,k=1,\ldots,n$, $j<k$.
\item $w_{ij} < w_{k\ell}$ \; if $i < k$, or, $i=k$ and $j<\ell$.
\end{itemize}
Also, each edge is of either of the following forms: 
\begin{enumerate}[(i)]
\item connecting $v_{i}$ and $w_{ij}$, or
\item connecting $v_{j}$ and $w_{ij}$.
\end{enumerate}
We assume that the type (i) edge has the signature $-$ and the other type has the signature $+$.

\begin{example}
The plumbing graphs for $n=2,3,4$ are described as follows.
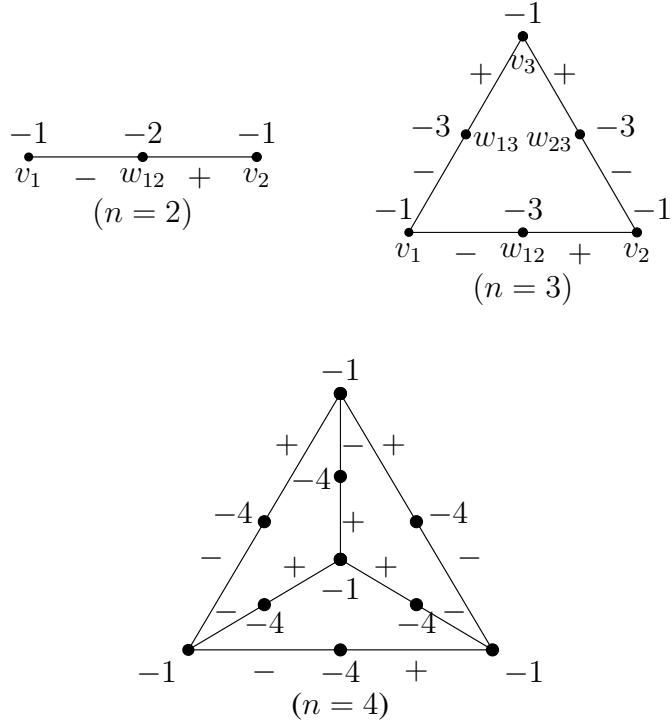
\[
\begin{tikzpicture}
\fill (0,0) circle (0.06);
\filldraw (0,0) --++(1.5,0) circle (0.06);
\filldraw (1.5,0) --++(1.5,0) circle (0.06);
\draw (0,0) node[above]{$-1$};
\draw (1.5,0) node[above]{$-2$};
\draw (3,0) node[above]{$-1$};

\draw (0,0) node[below]{$v_{1}$};
\draw (1.5,0) node[below]{$w_{12}$};
\draw (3,0) node[below]{$v_{2}$};

\draw (0.75,0) node[below]{$-$};
\draw (2.25,0) node[below]{$+$};

\draw (1.5,-0.75) node{$(n=2)$};

\coordinate (A) at (5,-1);

\fill (A) circle (0.06);
\filldraw (A) --++(1.5,0) circle (0.06);
\filldraw (A)++ (1.5,0) --++(1.5,0) circle (0.06);

\filldraw (A) --++(0.75,1.3) circle (0.06);
\filldraw (A) ++(0.75,1.3) --++(0.75,1.3) circle (0.06);

\filldraw (A) ++(3,0) --++(-0.75,1.3) circle (0.06);
\filldraw (A) ++(3,0) ++(-0.75,1.3) --++(-0.75,1.3) circle (0.06);

\draw (A) ++ (0,0)++(-0.2,0) node[above]{$-1$};
\draw (A) ++ (1.5,0) node[above]{$-3$};
\draw (A) ++ (3,0) ++(0.2,0) node[above]{$-1$};

\draw (A) ++ (0,0) node[below]{$v_{1}$};
\draw (A) ++ (1.5,0) node[below]{$w_{12}$};
\draw (A) ++ (3,0) node[below]{$v_{2}$};

\draw (A) ++ (0.75,0) node[below]{$-$};
\draw (A) ++ (2.25,0) node[below]{$+$};

\draw (A) ++ (1.5,2.6) node[above]{$-1$};
\draw (A) ++ (1.5,2.6) ++(0,-0.1)node[below]{$v_{3}$};

\draw (A) ++ (0.75,1.3)++(-0.05,-0.1) node[right]{$w_{13}$};
\draw (A) ++ (3,0) ++ (-0.75,1.3) ++(0.05,-0.1) node[left]{$w_{23}$};
\draw (A) ++ (0.75,1.3)++(-0.05,0.1) node[left]{$-3$};
\draw (A) ++ (3,0) ++ (-0.75,1.3) ++(0.05,0.1) node[right]{$-3$};

\draw (A) ++ (0.2,0.8) node{$-$};
\draw (A) ++ (0.2,0.8) ++ (0.75,1.3) node{$+$};

\draw (A) ++(3,0) ++ (-0.2,0.8) node{$-$};
\draw (A) ++(3,0) ++ (-0.2,0.8) ++ (-0.75,1.3) node{$+$};

\draw (A) ++ (1.5,-0.75) node{$(n=3)$};
\end{tikzpicture}
\]
\[
\begin{tikzpicture}

\fill (0,0) circle (0.08);
\filldraw (0,0) --++(2,0) circle (0.08);
\filldraw (2,0) --++(2,0) circle (0.08);

\filldraw (0,0) --++(1,1.7) circle (0.08);
\filldraw (1,1.7) --++(1,1.7) circle (0.08);

\filldraw (4,0) --++(-1,1.7) circle (0.08);
\filldraw (4,0) ++(-1,1.7) --++(-1,1.7) circle (0.08);

\fill (2,1.2) circle (0.08);

\filldraw (0,0) --++(1,0.6) circle (0.08);
\filldraw (1,0.6) --++(1,0.6) circle (0.08);

\filldraw (4,0) --++(-1,0.6) circle (0.08);
\filldraw (4,0) ++ (-1,0.6) --++(-1,0.6) circle (0.08);

\filldraw (2,3.4) --++(0,-1.1) circle (0.08);
\filldraw (2,3.4) ++(0,-1.1) --++(0,-1.1) circle (0.08);

\draw (1,0) node[below]{$-$};
\draw (3,0) node[below]{$+$};

\draw (0.3,1.5) node[below]{$-$};
\draw (1.3,3) node[below]{$+$};

\draw (4,0) ++ (-0.3,1.5) node[below]{$-$};
\draw (4,0) ++ (-1.3,3) node[below]{$+$};

\draw (2.17,2.7) node{$-$};
\draw (2.17,1.7) node{$+$};

\draw (0.5,0.5) node{$-$};
\draw (1.4,1.1) node{$+$};

\draw (4,0) ++ (-0.5,0.5) node{$-$};
\draw (4,0) ++ (-1.4,1.1) node{$+$};

\draw (0,0) node[below left]{$-1$};
\draw (4,0) node[below right]{$-1$};
\draw (2,3.4) node[above]{$-1$};
\draw (2,1.1)node[below]{$-1$};

\draw (2,0)node[below]{$-4$};
\draw (1,1.8) node[left]{$-4$};
\draw (3,1.8) node[right]{$-4$};
\draw (2.03,2.25) node[left]{$-4$};
\draw (1,0.65) node[below]{$-4$};
\draw (3,0.65) node[below]{$-4$};

\draw (2,-0.75) node{($n=4$)};
\end{tikzpicture}
\]

For each case, the weighted incidence graph $A_{n}$ is expressed as follows.
\[
A_{2} = 
\begin{pmatrix}
-1 & 0 & -1 \\
0 & -1 & 1 \\
-1 & 1 & -2 
\end{pmatrix} , \quad
A_{3} = 
\begin{pmatrix}
-1 & 0 & 0 & -1 & -1 & 0  \\
0 & -1 & 0 & 1 & 0 & -1 \\
0 & 0 & -1 & 0 & 1 & 1 \\
-1 & 1 & 0 & -3 & 0 & 0 \\
-1 & 0 & 1 & 0 & -3 & 0 \\
0 & -1 & 1 & 0 & 0 & -3 
\end{pmatrix} ,
\]
\[
A_{4} = 
\begin{pmatrix}
-1 & 0 & 0 & 0 & -1 & -1 & -1 & 0 & 0 & 0  \\
0 & -1 & 0 & 0 & 1 & 0 & 0 & -1 & -1 & 0  \\
0 & 0 & -1 & 0 & 0 & 1 & 0 & 1 & 0 & -1  \\
0 & 0 & 0 & -1 & 0 & 0 & 1 & 0 & 1 & 1  \\
-1 & 1 & 0 & 0 & -4 & 0 & 0 & 0 & 0 & 0  \\
-1 & 0 & 1 & 0 & 0 & -4 & 0 & 0 & 0 & 0  \\
-1 & 0 & 0 & 1 & 0 & 0 & -4 & 0 & 0 & 0  \\
0 & -1 & 1 & 0 & 0 & 0 & 0 & -4 & 0 & 0  \\
0 & -1 & 0 & 1 & 0 & 0 & 0 & 0 & -4 & 0  \\
0 & 0 & -1 & 1 & 0 & 0 & 0 & 0 & 0 & -4  
\end{pmatrix} .
\]
\end{example}

In this section, we will prove the following.
\begin{theorem}\label{thm:main}
\[
H_{1} (\partial \overline{F} , \bZ) \cong \bZ^{n(n-1)/2} \oplus \bZ_{n}^{(n-2)(n-3)/2}.
\]
\end{theorem}
From Theorem \ref{thm:homology}, it suffices to compute the Smith normal form of the weighted incidence matrix.
For matrices $D_{1} , \ldots, D_{\ell}$, we denote their direct sum by $D_{1} \oplus \cdots \oplus D_{\ell}$, that is,
\[
D_{1} \oplus \cdots \oplus D_{\ell} = 
\begin{pmatrix}
D_{1} & O & \cdots & O \\
O & D_{2} & \cdots & O \\
\vdots & \vdots & \ddots & \vdots \\
O & O & \cdots & D_{\ell} 
\end{pmatrix},
\]

where $O$ is the zero matrix.
Under the fixed order of vertices as above, the weighted incidence matrix $A_{n}$ is expressed as:
\[
A_{n} = 
\begin{pmatrix}
-E_{n} & -^{t}X_{n} \\
-X_{n} & -nE_{n(n-1)/2}
\end{pmatrix},
\]

where $E_{k}$ is the identity matrix of size $k$ and $X_{n}$ is a matrix of size $(n(n-1)/2,n)$ defined inductively as follows:
\[
X_n = 
\begin{pmatrix}
\bf{1} & -E_{n-1} \\
\bf{0} & X_{n-1} 
\end{pmatrix} \, \mbox{($n>2$)}, \quad
X_{2} = 
\begin{pmatrix}
1 & -1 
\end{pmatrix}.
\]
Here,  $\bf{0}$ is the zero vector and $\bf{1}$ is the column vector with all entries equal to $1$.
For example, $X_3$, $X_4$ are expressed as follows:
\[
X_{3} = 
\begin{pmatrix}
1 & -1 & 0 \\
1 & 0 & -1 \\
0 & 1 & -1 
\end{pmatrix} , \quad
X_{4} = 
\begin{pmatrix}
1 & -1 & 0 & 0  \\
1 & 0 & -1 & 0  \\
1 & 0 & 0 & -1 \\
0 & 1 & -1 & 0 \\
0 & 1 & 0 & -1 \\
0 & 0 & 1 & -1 
\end{pmatrix} .
\]

We set $B_{n} = nE_{n(n-1)/2} - X_{n}\, ^{t}X_{n}$. 
Since we have 
\begin{eqnarray*}
& &
-A_{n}
\begin{pmatrix}
E_{n} & -^{t} X_{n} \\
O & E_{n(n-1)/2}
\end{pmatrix} \\
&=& 
\begin{pmatrix}
E_{n} & {}^{t} X_{n} \\
X_{n} & nE_{n(n-1)/2}
\end{pmatrix}
\begin{pmatrix}
E_{n} & -^{t} X_{n} \\
O & E_{n(n-1)/2}
\end{pmatrix} \\
&=&
\begin{pmatrix}
E_{n} & O \\
X_{n} & B_{n}
\end{pmatrix}, 
\end{eqnarray*} 
the matrix $A_{n}$ changes to $E_{n} \oplus B_{n}$ by an elementary transformation. 

\begin{lemma}\label{lem:keylem}
Let $J_{k}$ be a square matrix of size $k$ with all entries equal to $1$.
\begin{enumerate}[(1)]
\item The matrix $X_{n}\, ^{t}X_{n}$ has the following block matrix decomposition:
\[
X_{n}\, ^{t}X_{n} = 
\begin{pmatrix}
J_{n-1} + E_{n-1} & -^{t} X_{n-1} \\
-X_{n-1} & X_{n-1} \, ^{t}X_{n-1}
\end{pmatrix}.
\]
\item $^{t}X_{n} X_{n} = nE_{n} - J_{n}$.
\item The matrix $B_{n}$ has the following block matrix decomposition:
\[
B_{n} = 
\begin{pmatrix}
^{t}X_{n-1}X_{n-1} & ^{t}X_{n-1} \\
X_{n-1} & nE_{(n-1)(n-2)/2} -X_{n-1} \,^{t}X_{n-1}
\end{pmatrix}.
\]
\item $(nE_{n(n-1)/2} -X_{n}\, ^{t}X_{n}) X_{n} = O$. Particularly, we have 
\[
(nE_{(n-1)(n-2)/2} - X_{n-1} \,^{t}X_{n-1})X_{n-1} = X_{n-1}.
\]
\end{enumerate}
\end{lemma}

\proof
We will omit the size of the identity matrix in proofs except in special cases.
\begin{enumerate}[(1)]
\item It follows from direct computation:
\begin{eqnarray*}
X_{n} \, ^{t} X_{n} &=& 
\begin{pmatrix}
{\bf 1} & -E \\
{\bf 0} & X_{n-1} 
\end{pmatrix}
\begin{pmatrix}
{}^{t} {\bf 1} & {}^{t} {\bf 0} \\
-E & ^{t} X_{n-1} 
\end{pmatrix}
\\
&=&
\begin{pmatrix}
J_{n-1} + E & -^{t} X_{n-1} \\
-X_{n-1} & X_{n-1} \, ^{t} X_{n-1}
\end{pmatrix}.
\end{eqnarray*}

\item We prove by induction. Clearly, it holds when $n=2$. When $n>2$, we have
\begin{eqnarray*}
^{t}X_{n} X_{n} &=&
\begin{pmatrix}
{}^{t} {\bf 1} & {}^{t} {\bf 0} \\
-E & ^{t} X_{n-1}
\end{pmatrix}
\begin{pmatrix}
{\bf 1} & -E \\
{\bf 0} & X_{n-1}
\end{pmatrix} \\
&=&
\begin{pmatrix}
n-1 & -{}^{t} {\bf 1}\\
- {\bf 1} & E + \, ^{t}X_{n-1} X_{n-1}
\end{pmatrix} \\
&=& 
\begin{pmatrix}
n-1 & -{}^{t} {\bf 1}\\
- {\bf 1} & E + (n-1)E - J_{n-1}
\end{pmatrix} \\
&=& nE - J_{n}.
\end{eqnarray*}

\item Using (1) and (2), we have 
\begin{eqnarray*}
B_{n} &=& nE -X_{n} \, ^{t}X_{n} \\
&=& nE -
\begin{pmatrix}
J_{n-1} + E & -^{t} X_{n-1} \\
-X_{n-1} & X_{n-1} \, ^{t}X_{n-1}
\end{pmatrix} \\
&=&
\begin{pmatrix}
(n-1)E - J_{n-1} & ^{t} X_{n-1} \\
X_{n-1} & nE - X_{n-1} \, ^{t}X_{n-1}
\end{pmatrix} \\
&=&
\begin{pmatrix}
^{t}X_{n-1} X_{n-1} & ^{t} X_{n-1} \\
X_{n-1} & nE - X_{n-1} \, ^{t}X_{n-1}
\end{pmatrix}.
\end{eqnarray*}

\item
It follows from the following computation.
Note that the sum of each row of $X_{n}$ is equal to $0$.
\begin{eqnarray*}
&& (nE - X_{n} \, ^{t}X_{n})X_{n} \\
&= &X_{n} (nE - \, ^{t} X_{n} X_{n})  \\
&= &X_{n} (nE - (nE -J_{n}) ) \\
&=& X_{n} \cdot J_{n} \\
&=& O.
\end{eqnarray*}
\end{enumerate}
\endproof

\begin{theorem}\label{thm:key}
The Smith normal form of $A_{n}$ is 
\[
E_{2n-2} \oplus nE_{(n-2)(n-3)/2} \oplus O.
\]
\end{theorem}

\proof
It suffices to compute the Smith normal form of $B_{n}$.
From Lemma \ref{lem:keylem}, we have the following matrix by performing an elementary row transformation to $B_{n}$:
\begin{eqnarray*} 
&&
\begin{pmatrix}
E & - \, ^{t}X_{n-1} \\
O & E
\end{pmatrix}
B_{n} \\
&=& 
\begin{pmatrix}
E & - \, ^{t}X_{n-1} \\
O & E
\end{pmatrix}
\begin{pmatrix}
^{t}X_{n-1} X_{n-1} & ^{t} X_{n-1} \\
X_{n-1} & nE - X_{n-1} \, ^{t}X_{n-1}
\end{pmatrix}\\
&=& 
\begin{pmatrix}
O & ^{t}X_{n-1} - \,^{t}X_{n-1}(nE - X_{n-1} \, ^{t}X_{n-1} )\\
X_{n-1} & nE -X_{n-1} \, ^{t}X_{n-1}
\end{pmatrix} \\
&=&
\begin{pmatrix}
O & O \\
X_{n-1} & nE - X_{n-1} \,^{t}X_{n-1} 
\end{pmatrix}.
\end{eqnarray*}
Next, by reordering the rows and performing an elementary column transformation, we have
\begin{eqnarray*}
&&
\begin{pmatrix}
X_{n-1} & nE - X_{n-1} \,^{t}X_{n-1} \\
O & O 
\end{pmatrix} 
\begin{pmatrix}
E & ^{t} X_{n-1} \\
O & E
\end{pmatrix}
\\
&=& 
\begin{pmatrix}
X_{n-1} & nE \\
O & O \\
\end{pmatrix} \\
&=&
\begin{pmatrix}
{\bf 1} & -E_{n-2} & nE_{n-2} & O \\
{\bf 0} & X_{n-2} & O & nE_{(n-2)(n-3)/2}  \\
{\bf 0} & O & O & O
\end{pmatrix}.
\end{eqnarray*}
Finally, by row and column reduction, the Smith normal form of $B_{n}$ is expressed as 
\[
\begin{pmatrix}
E_{n-2} & O & O \\
O & nE_{(n-2)(n-3)/2} & O \\
O & O & O
\end{pmatrix}.
\]
The theorem follows from $A_{n} = E_{n} \oplus B_{n}$.
\endproof

\textbf{Proof of Theorem \ref{thm:main}.}
We apply Theorem \ref{thm:homology}. 
The Betti number of the plumbing graph is equal to the one of the complete graph, thus $b_{1} (\Gamma) = (n-1)(n-2)/2$.
By Theorem \ref{thm:key}, corank of $A_{n}$ is $n-1$ and the torsion part is $\bZ_{n}^{(n-2)(n-3)/2}$. Hence we have
\[
H_{1} (\partial F, \bZ) = \bZ^{n(n-1)/2} \oplus \bZ_n^{(n-2)(n-3)/2}.
\]
\qed

\section{Plumbing graphs of the Milnor fiber boundary for arrangements}\label{sec:algo}
In this section, we will give the algorithm to obtain the plumbing graph of the Milnor fiber boundary for hyperplane arrangements in $\bC^3$, following the book \cite{nem-szi}.
Roughly speaking, we need the following three steps to obtain the plumbing graph in general. 
Let $f: (\bC^3, 0) \rightarrow ( \bC, 0)$ be a defining function germ of a non-isolated surface singularity.
\begin{enumerate}[(i)]
\item Find a generic holomorphic function germ $g: (\bC^3,0) \rightarrow (\bC, 0)$
such that the pair $(f,g)$ defines an isolated complete intersection singularity pair (shortly, ICIS pair) (Chapter 3 in \cite{nem-szi}).
\item Compute the dual graph $\Gamma_{\cC}$ of a curve configuration determined by the ICIS pair $(f,g)$ (Chapter 6 in \cite{nem-szi}).
\item Construct the desired plumbing graph from $\Gamma_{\cC}$ with a procedure called the ``main algorithm" (Chapter 10 in \cite{nem-szi}).
\end{enumerate}

However, the complete explanation needs a lot of preparation. 
Thus, we accept the structure of the graph $\Gamma_{\cC}$ in the case of arrangements and the main algorithm to obtain the plumbing graph as fact.
We then focus on specific calculations.

\subsection{The curve configuration graph}

The graph $\Gamma_{\cC}$ is determined by the ICIS pair $(f,g)$.
We do not explain the general definition and construction of $\Gamma_{\cC}$ from a given non-isolated surface singularity and ICIS pair $(f,g)$ (see Chapter 6 in \cite{nem-szi} for details).
Instead, we will explain some properties of $\Gamma_{\cC}$ and the construction for a given hyperplane arrangement in $\bC^3$.

The graph $\Gamma_{\cC}$ has a special kind of vertex called \textit{arrowheads}. An arrowhead is a degree one vertex.
The edge connecting an arrowhead and another vertex is called an \textit{arrow}.
The vertex set is decomposed into arrowhead vertices and non-arrowhead vertices.
The graph $\Gamma_{\cC}$ has the following decoration:
\begin{itemize}
\item Each vertex is decorated by a triple $(m;n,\nu)$ of nonnegative integers.
\item Each edge is decorated by the number $1$ or $2$.
\end{itemize}

Moreover, it is known that the graph $\Gamma_{\cC}$ has only edges and vertices with the following form:
\begin{itemize}
\item Type $1$-edge:
	\[
	\begin{tikzpicture}
	\draw (0,0)--++(2,0);
	\fill (0,0) circle (0.06);
	\fill (2,0) circle (0.06);
	\draw (0,0) node[above]{$(m;n,\nu)$};
	\draw (2,0) node[above]{$(m;\ell,\lambda)$};
	\draw (1,0) node[below]{$1$};
	\end{tikzpicture}
	\]
\item Type $2$-edge:
	\[
	\begin{tikzpicture}
	\draw (0,0)--++(2,0);
	\fill (0,0) circle (0.06);
	\fill (2,0) circle (0.06);
	\draw (0,0) node[above]{$(m;n,\nu)$};
	\draw (2,0) node[above]{$(m';n,\nu)$};
	\draw (1,0) node[below]{$2$};
	\end{tikzpicture}
	\]
\end{itemize}
In other words, if two vertices are connected by a $1$-edge, then the first entries coincide. If two vertices are connected by a $2$-edge, then the second and third entries coincide.
An arrowhead vertex is always decorated by $(1;0,1)$ and an arrow edge is always decorated by $1$:
	\[
	\begin{tikzpicture}
	\draw[->] (0,0)--++(2,0);
	\fill (0,0) circle (0.06);
	\draw (0,0) node[above]{$(1;n,\nu)$};
	\draw (2,0) node[above]{$(1;0,1)$};
	\draw (1,0) node[below]{$1$};
	\end{tikzpicture}
	\]
In the case of a hyperplane arrangement, we can construct $\Gamma_{\cC}$ as the following theorem (see Proposition 8.1.1 and Subsection 8.2 in \cite{nem-szi}). 
Let $\cA$ be a hyperplane arrangement in $\bC^3$.
Recall that $\Gamma(\overline{\cA})$ is the incidence graph of the corresponding projective line arrangement.  
Let $n= |\cA|$. For each intersection point $p_{j}$, let $m_{j} = \# \{\overline{H} \in \overline{\cA} \mid \overline{H} \ni p_{j}\}$.

\begin{theorem}\label{thm:dual}
The graph $\Gamma_{\cC}$ for $\cA$ is obtained from $\Gamma(\overline{\cA})$ by adding an arrow to each line vertex.
Each line vertex $v_{i}$ is decorated by $(1;n,1)$ and each point vertex $w_{j}$ is decorated by $(m_{j};n,1)$. All non-arrow edges are decorated by $2$.
\end{theorem}

\begin{example} \textbf{A pencil arrangement.}
Suppose that $\cA$ is defined by a polynomial $Q = x^n + y^n$.
Then, all lines in $\overline{\cA}$ have a common intersection point $p=[0:0:1]$.
The graph $\Gamma_{\cC}$ is described as:
	
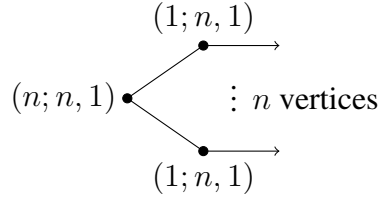
\[
	\begin{tikzpicture}
	\filldraw (0,0) circle (0.06) node[left]{$(n;n,1)$};
	\draw[->] (0,0) --++(1,0.7) --++(1,0);
	\draw[->] (0,0) --++(1,-0.7) --++(1,0);
	\filldraw (1,0.7) circle (0.06) node[above]{$(1;n,1)$};
	\filldraw (1,-0.7) circle (0.06) node[below]{$(1;n,1)$};
	\draw (1.4,0.1) node{$\vdots$};
	\draw (1.5,0) node[right]{$n$ vertices};
	\end{tikzpicture}
	\]
\end{example}

\begin{example}
\textbf{A near pencil arrangement.}
Suppose that $\cA$ is defined by a polynomial $Q = z (x^{n-1} - y^{n-1})$.
Let $H_{0} = \{z=0\}$ and $H_{j} = \{x - e^{2 \pi j \sqrt{-1}/(n-1)} y = 0\}$.
Let $p_{0} = [0:0:1] , p_{j}=[e^{2 \pi j \sqrt{-1}/(n-1)}:1:0]$ be the intersection points in $\overline{\cA}$.
We have $p_{0} = H_{1} \cap \cdots \cap H_{n-1}$ and $p_{j} =H_{0} \cap H_{j}$ for $j=1, \ldots, n-1$.
The graph $\Gamma_{\cC}$ is described as:
	
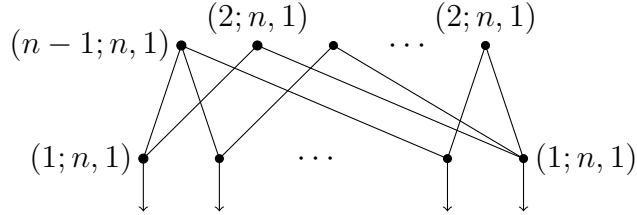
\[
	\begin{tikzpicture}
	\filldraw (0,0) circle (0.06) node[left]{$(1;n,1)$};
	\fill (1,0) circle (0.06) ;
	\fill (4,0) circle (0.06) ;
	\fill (5,0) circle (0.06) node[right]{$(1;n,1)$};
	\filldraw (0.5,1.5) circle (0.06) node[left]{$(n-1;n,1)$};	
	\filldraw (1.5,1.5) circle (0.06) node[above]{$(2;n,1)$};
	\fill (2.5,1.5) circle (0.06) ;		
	\fill (4.5,1.5) circle (0.06) node[above]{$(2;n,1)$};
	\draw[->] (0,0) --++(0,-0.7);	
	\draw[->] (1,0) --++(0,-0.7);
	\draw[->] (4,0) --++(0,-0.7);
	\draw[->] (5,0) --++(0,-0.7);
	\draw (0,0) --(0.5,1.5);		
	\draw (0,0) --(1.5,1.5);
	\draw (1,0) --(0.5,1.5);	
	\draw (1,0) --(2.5,1.5);
	\draw (4,0) --(0.5,1.5);	
	\draw (4,0) --(4.5,1.5);
	\draw (5,0) --(1.5,1.5);
	\draw (5,0) --(2.5,1.5);
	\draw (5,0) --(4.5,1.5);
	\draw (2.3,0) node{$\cdots$};
	\draw (3.5,1.5) node{$\cdots$};		
	\end{tikzpicture}
	\]
\end{example}

\begin{example}
\textbf{A generic arrangement.}
Let $\cA$ be a generic arrangement. 
At that time, every intersection point in $\overline{\cA}$ is the intersection of exactly two lines.
The graph $\Gamma_{\cC}$ is constructed from a complete graph of $n$ vertices (line vertex)  by inserting one vertex (point vertex) to each edge and adding an arrow from each line vertex.
For example, the graph $\Gamma_{\cC}$ for $n=4$ is described as follows:
	
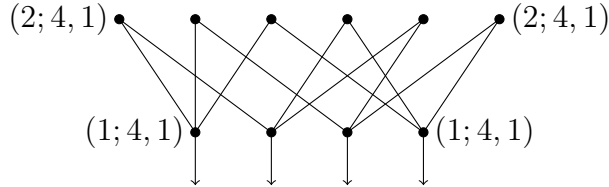
\[
	\begin{tikzpicture}
	\filldraw (0,0) circle (0.06) node[left]{$(1;4,1)$};
	\filldraw (1,0) circle (0.06);
	\filldraw (2,0) circle (0.06);
	\filldraw (3,0) circle (0.06) node[right]{$(1;4,1)$};
	\filldraw (-1,1.5) circle (0.06) node[left]{$(2;4,1)$};
	\filldraw (0,1.5) circle (0.06) ;
	\filldraw (1,1.5) circle (0.06) ;
	\filldraw (2,1.5) circle (0.06) ;
	\filldraw (3,1.5) circle (0.06) ;	
	\filldraw (4,1.5) circle (0.06) node[right]{$(2;4,1)$};	
	\draw[->] (0,0) --++(0,-0.7);
	\draw[->] (1,0) --++(0,-0.7);
	\draw[->] (2,0) --++(0,-0.7);
	\draw[->] (3,0) --++(0,-0.7);
	\draw (0,0) --(-1,1.5);			
	\draw (0,0) --(0,1.5);
	\draw (0,0) --(1,1.5);	
	\draw (1,0) --(-1,1.5);		
	\draw (1,0) --(2,1.5);	
	\draw (1,0) --(3,1.5);
	\draw (2,0) --(0,1.5);	
	\draw (2,0) --(3,1.5);
	\draw (2,0) --(4,1.5);	
	\draw (3,0) --(1,1.5);
	\draw (3,0) --(2,1.5);
	\draw (3,0) --(4,1.5);	
	\end{tikzpicture}
	\]
\end{example}

\subsection{The string graph $Str^{\pm} (a,b;c \mid i,j;k)$}
Before going to the main algorithm, we need one more preparation. 
The main algorithm requires a special form of a string graph $Str^{\pm}(a,b;c \mid i,j;k)$ for nonnegative integers with $\gcd(a,b,c)=1$.
We will give the construction (see Section 4.3.9 in \cite{nem-szi}).
In the case of arrangements, all the possible graphs to appear satisfy the condition $(i,j,k) = (0,0,1)$. Thus, we assume this condition. 
The notation $(a,c)$ stands for $\gcd(a,c)$.

First, we consider the unique integers $0 \leq \lambda < c /(a,c)$ and $m_{1} \in \bZ_{>0}$ with the following equation:
\begin{equation} \label{eq:lambda}
b + \lambda \cdot \frac{a}{(a,c)} = m_1 \cdot \frac{c}{(a,c)}
\end{equation}
If $\lambda \neq 0$, we next consider the continued fraction:
\[
\frac{c/(a,c)}{\lambda} = k_{1} - \cfrac{1}{k_2 - \cfrac{1}{\ddots - \cfrac{1}{k_{s}}}} , \quad k_{1} ,\ldots, k_{s} \geq 2.
\]
Using the above integers, we define $m_{2}, \ldots, m_{s}$ as:
\begin{eqnarray*}
m_{2} &=& k_{1} m_{1} - \frac{a}{(a,c)}, \\
m_{i+1} &=& k_{i}m_{i} - m_{i-1}, \quad \mbox{for $i \geq 2$.}
\end{eqnarray*}
Then, the graph $Str(a,b;c \mid 0,0;1)$ is described as follows:
	
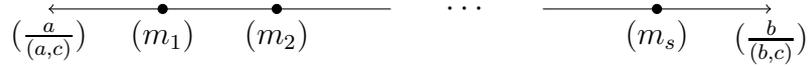
\[
	\begin{tikzpicture}
	\filldraw (0,0) circle (0.06) node[below]{$(m_{1})$};
	\filldraw (1.5,0) circle (0.06)node[below]{$(m_{2})$};
	\draw[->] (0,0) --++ (-1.5,0)node[below]{$(\frac{a}{(a,c)})$};
	\draw (0,0) --++(3,0);
	\draw (4,0) node{$\cdots$};
	\draw[->] (5,0) --++(3,0) node[below]{$(\frac{b}{(b,c)})$};
	\filldraw (6.5,0) circle (0.06) node[below]{$(m_{s})$};
	\end{tikzpicture}
	\]
the integer $(m_{i})$ is the multiplicity weight of each vertex.The graph $Str^{-} (a,b;c \mid 0,0;1)$ is the same string with all edges decorated by the signature $-$.

If $\lambda = 0$, then the graph $Str^{\pm}(a,b;c \mid 0,0;1)$ is a double arrow having no non-arrowhead vertex, with the edge decorated by $\pm$.

We show some examples. Let $n \geq 2$ and $k \geq 1$ be integers. 
\begin{example}\label{exam:str:(n,n)}
$Str(1,n;n \mid 0,0;1)$.

The equation (\ref{eq:lambda}) becomes
\[
n + \lambda \cdot \frac{1}{(1,n)} = m_{1} \cdot \frac{n}{(1,n)}.
\]
Since we have $\lambda = 0$, the string graph $Str(1,n;n\mid 0,0;1)$ is as follows:
	
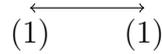
\[
	\begin{tikzpicture}		
	\draw[<->] (0,0) --++ (-1.5,0)node[below]{$(1)$};
	\draw (0,0) node[below]{$(1)$};
	\end{tikzpicture}
	\]
\end{example}

\begin{example}\label{exam:str:(n-1,n)}
$Str(1,n-1;n\mid 0,0;1)$.

The equation (\ref{eq:lambda}) becomes
\[
n-1 + \lambda \cdot \frac{1}{(1,n)} = m_{1} \cdot \frac{n}{(1,n)}.
\]
Thus we have $\lambda = 1$, $m_{1} = 1$. The continued fraction is 
\[
\frac{n/(1,n)}{1} = n.
\]
Therefore, the string graph $Str(1,n-1;n\mid 0,0;1)$ is as follows:
	
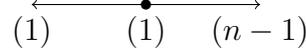
\[
	\begin{tikzpicture}
	\draw[<->] (0,0) --++(3,0) node[below]{$(n-1)$};
	\draw (0,0) node[below]{$(1)$};
	\filldraw (1.5,0) circle (0.06) node[below]{$(1)$};
	\end{tikzpicture}
	\]
\end{example}

\begin{example}\label{exam:str:(2,2k)}
$Str(1,2;2k \mid 0,0;1)$.

The solution of equation (\ref{eq:lambda}) is $\lambda = 2k-2$, $m_{1} = 1$.
The continued fraction is 
\[
\frac{2k}{2k-2} = 2 - \cfrac{1}{2 - \cfrac{1}{\ddots - \cfrac{1}{2}}}, \quad \mbox{($(k-1)$-times)}
\]
The multiplicity $m_{i} = 1$ for all $i$. Thus, the string graph $Str(1,2;2k\mid 0,0;1)$ is as follows:
	\[
	\begin{tikzpicture}
	\filldraw (0,0) circle (0.06) node[below]{$(1)$};
	\filldraw (1.5,0) circle (0.06)node[below]{$(1)$};
	\draw[->] (0,0) --++ (-1.5,0)node[below]{$(1)$};
	\draw (0,0) --++(3,0);
	\draw (4,0) node{$\cdots$};
	\draw[->] (5,0) --++(3,0) node[below]{$(1)$};
	\filldraw (6.5,0) circle (0.06) node[below]{$(1)$};
	\draw (3.5,-0.7) node[below]{($k-1$ non-arrowhead vertices)};
	\end{tikzpicture}
	\]
\end{example}

\begin{example}\label{exam:str:(2,2k+1)}
$Str(1,2;2k+1 \mid 0,0;1)$.

The solution of equation (\ref{eq:lambda}) is $\lambda = 2k-1$, $m_{1} = 1$.
The continued fraction is 
\[
\frac{2k+1}{2k-1} = 2 - \cfrac{1}{\ddots - \cfrac{1}{2- \cfrac{1}{3}}}, \quad 
\]
where $2$ appears $(k-1)$ times. All the multiplicities are all equal to $1$. Thus, the string graph $Str(1,2;2k+1 \mid 0,0;1)$ is as follows:
	
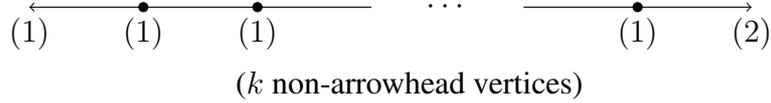
\[
	\begin{tikzpicture}
	\filldraw (0,0) circle (0.06) node[below]{$(1)$};
	\filldraw (1.5,0) circle (0.06)node[below]{$(1)$};
	\draw[->] (0,0) --++ (-1.5,0)node[below]{$(1)$};
	\draw (0,0) --++(3,0);
	\draw (4,0) node{$\cdots$};
	\draw[->] (5,0) --++(3,0) node[below]{$(2)$};
	\filldraw (6.5,0) circle (0.06) node[below]{$(1)$};
	\draw (3.5,-0.7) node[below]{($k$ non-arrowhead vertices)};
	\end{tikzpicture}
	\]
\end{example}

\begin{remark}
In the case of arrangements, only the type $Str^{-}(1,m;n \mid 0,0;1)$ appears ($n \geq m \geq 2$) as we will see later. Thus, we always have $\lambda=n-m$ and $m_{1} = 1$.
\end{remark}


\subsection{The main algorithm specialized in arrangements}
We now give how to construct the plumbing graph of the Milnor fiber boundary for a hyperplane arrangement from the graph $\Gamma_{\cC}$.
The plumbing graph is constructed by modifying $\Gamma_{\cC}$ and giving decorations as follows (this is what is called the \textit{main algorithm}):

Vertices:
\begin{enumerate}[(i)]
\item 
For a line vertex $v_{i}$, the genus is $[0]$ and multiplicity is $(1)$. 

\item 
For a point vertex $w_{j}$, the genus is $[(m_{j} -2) (\gcd(m_{j},n)-1)/2]$ and the multiplicity is $(m_j / \gcd(m_{j},n))$.

\item 
For an arrowhead, the genus is $[0]$ and the multiplicity is $(1)$.
\end{enumerate}

Edges:
\begin{enumerate}[(i)]
\item 
A non-arrow edge has the form:
	
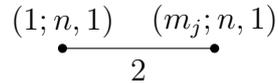
\[
	\begin{tikzpicture}
	\draw (0,0)--++(2,0);
	\fill (0,0) circle (0.06);
	\fill (2,0) circle (0.06);
	\draw (0,0) node[above]{$(1;n,1)$};
	\draw (2,0) node[above]{$(m_{j};n,1)$};
	\draw (1,0) node[below]{$2$};
	\end{tikzpicture}
	\]
We insert the string graph of type $Str^{-}(1,m_{j};n \mid 0,0;1)$ to the edge.

\item 
An arrow edge becomes an edge with the decoration $+$. Since it connects an arrowhead and a vertex, it is still an arrow. 
\end{enumerate}

The Euler number of a vertex is computed by using the formula (\ref{eq:mult}) in Section \ref{sec:graph}.
The above output is the result of the specific computation written in 103--105 pages of \cite{nem-szi} for the case of hyperplane arrangements.
Under this setting, we have the following theorem.
\begin{theorem} (Theorem 12.10 in \cite{nem-szi})
The plumbing graph of the Milnor fiber boundary is obtained by the above procedures and by removing arrowheads.
\end{theorem}

\begin{remark}
In general setting, we need the notion of covering graphs. 
However, we do not need it in the case of arrangements (we always have the covering degree $\mathfrak{n}_{w}, \mathfrak{n}_{e} = 1$ in terms of \cite{nem-szi}).
Thus we did not explain about it.
\end{remark}

We give some examples of the computation of the main algorithm.
\begin{example}
Suppose that $\cA$ is a pencil arrangement of $n$ hyperplanes. Recall that the graph $\Gamma_{\cC}$ is described as:
	
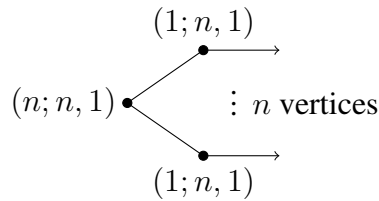
\[
	\begin{tikzpicture}
	\filldraw (0,0) circle (0.06) node[left]{$(n;n,1)$};
	\draw[->] (0,0) --++(1,0.7) --++(1,0);
	\draw[->] (0,0) --++(1,-0.7) --++(1,0);
	\filldraw (1,0.7) circle (0.06) node[above]{$(1;n,1)$};
	\filldraw (1,-0.7) circle (0.06) node[below]{$(1;n,1)$};
	\draw (1.4,0.1) node{$\vdots$};
	\draw (1.5,0) node[right]{$n$ vertices};
	\end{tikzpicture}
	\]
By applying the main algorithm, the string of type $Str(1,n;n \mid 0,0;1)$ is inserted to each vertex with trivial covering data, and we have
	
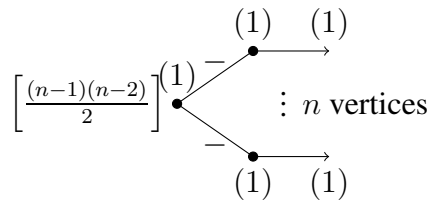
\[
	\begin{tikzpicture}
	\filldraw (0,0) circle (0.06) node[above]{$(1)$};
	\draw (-0,0) node[left]{$\left [ \frac{(n-1)(n-2)}{2} \right ]$};
	\draw[->] (0,0) --++(1,0.7) --++(1,0) node[above]{$(1)$};
	\draw[->] (0,0) --++(1,-0.7) --++(1,0) node[below]{$(1)$};
	\filldraw (1,0.7) circle (0.06) node[above]{$(1)$};
	\filldraw (1,-0.7) circle (0.06) node[below]{$(1)$};
	\draw (1.4,0.1) node{$\vdots$};
	\draw (1.5,0) node[right]{$n$ vertices};
	\draw (0.5,0.55) node{$-$};
	\draw (0.5,-0.55) node{$-$};
	\end{tikzpicture}
	\]
Next, by applying the formula (\ref{eq:mult}) and removing arrows, we have
	
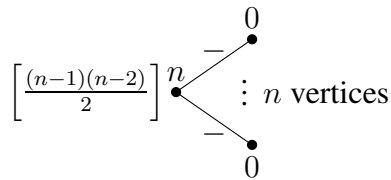
\[
	\begin{tikzpicture}
	\filldraw (0,0) circle (0.06) node[above]{$n$};
	\draw (-0,0) node[left]{$\left [ \frac{(n-1)(n-2)}{2} \right ]$};
	\draw (0,0) --++(1,0.7) ;
	\draw (0,0) --++(1,-0.7) ;
	\filldraw (1,0.7) circle (0.06) node[above]{$0$};
	\filldraw (1,-0.7) circle (0.06) node[below]{$0$};
	\draw (0.9,0.1) node{$\vdots$};
	\draw (1,0) node[right]{$n$ vertices};
	\draw (0.5,0.55) node{$-$};
	\draw (0.5,-0.55) node{$-$};
	\end{tikzpicture}
	\]
Finally, by plumbing calculus (splitting), we obtain
	\[
	\begin{tikzpicture}
	\draw (-0.3,0) node[left]{$(n-1)^2$ copies of};
	\filldraw (0,0) circle (0.06) node[above]{$0$};	
	\end{tikzpicture}
	\]
Therefore, $\partial \overline{F}$ is diffeomorphic to $\sharp_{(n-1)^2} S^1 \times S^2$.
\end{example}

\begin{example}
Suppose that $\cA$ is a near pencil arrangement with $n$ hyperplanes. 
The graph $\Gamma_{\cC}$ is described as
	
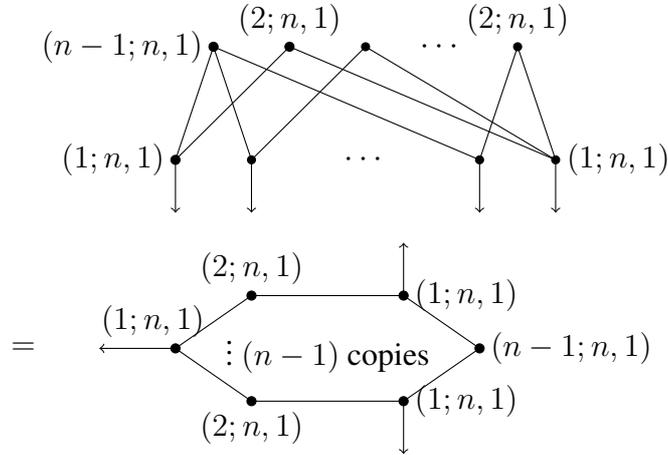
\[
	\begin{tikzpicture}
	\filldraw (0,0) circle (0.06) node[left]{$(1;n,1)$};
	\fill (1,0) circle (0.06) ;
	\fill (4,0) circle (0.06) ;
	\fill (5,0) circle (0.06) node[right]{$(1;n,1)$};
	\filldraw (0.5,1.5) circle (0.06) node[left]{$(n-1;n,1)$};	
	\filldraw (1.5,1.5) circle (0.06) node[above]{$(2;n,1)$};
	\fill (2.5,1.5) circle (0.06) ;		
	\fill (4.5,1.5) circle (0.06) node[above]{$(2;n,1)$};
	\draw[->] (0,0) --++(0,-0.7);	
	\draw[->] (1,0) --++(0,-0.7);
	\draw[->] (4,0) --++(0,-0.7);
	\draw[->] (5,0) --++(0,-0.7);
	\draw (0,0) --(0.5,1.5);		
	\draw (0,0) --(1.5,1.5);
	\draw (1,0) --(0.5,1.5);	
	\draw (1,0) --(2.5,1.5);
	\draw (4,0) --(0.5,1.5);	
	\draw (4,0) --(4.5,1.5);
	\draw (5,0) --(1.5,1.5);
	\draw (5,0) --(2.5,1.5);
	\draw (5,0) --(4.5,1.5);
	\draw (2.5,0) node{$\cdots$};
	\draw (3.5,1.5) node{$\cdots$};
	
	\draw (-2,-2.5)  node{$=$};
	
	\filldraw (0,-2.5) circle (0.06) ;
	\draw (-0.3,-2.5) node[above]{$(1;n,1)$};
	\draw [->] (0,-2.5) --++(-1,0);
	\draw (0,-2.5) --++(1,0.7);
	\draw (0,-2.5) --++(1,-0.7);
	\filldraw (1,-1.8) circle (0.06) node[above]{$(2;n,1)$};
	\draw (1,-1.8) --++(2,0);
	\draw[->] (1,-1.8) ++(2,0) --++(0,0.7);
	\filldraw (1,-1.8) ++(2,0) circle (0.06) node[right]{$(1;n,1)$};
	\filldraw (1,-3.2) circle (0.06) node[below]{$(2;n,1)$};
	\draw (1,-3.2) --++(2,0);
	\draw[->] (1,-3.2) ++(2,0) --++(0,-0.7);	
	\filldraw (1,-3.2) ++(2,0) circle (0.06) node[right]{$(1;n,1)$};
	\draw (1,-1.8) ++(2,0) --++(1,-0.7);
	\draw (1,-3.2) ++(2,0) --++(1,0.7);
	\filldraw (1,-3.2) ++(2,0) ++(1,0.7) circle (0.06) node[right]{$(n-1;n,1)$};
	\draw (2,-2.5) node{$\vdots$ $(n-1)$ copies};
	\end{tikzpicture}
	\]

By applying the main algorithm, strings of type $Str^{-}(1,n-1;n \mid 0,0;1)$ and $Str^{-}(1,2;2k \mid 0,0;1)$ ($n=2k$) or $Str^{-}(1,2;2k+1 \mid 0,0;1)$ ($n=2k+1$) are inserted with trivial covering data, and we have
	
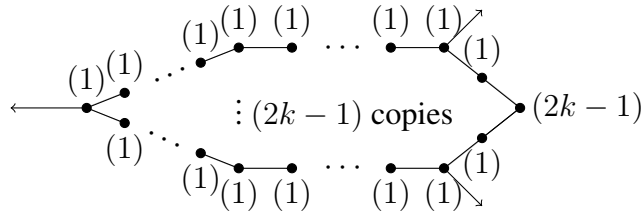
\[
	\begin{tikzpicture}
	\draw[->] (0,0) --++(-1,0);
	\filldraw (0,0) circle (0.06) node[above]{$(1)$};
	\filldraw (0,0) --++(0.5,0.2) circle (0.06) node[above]{$(1)$};
	\draw (1,0.5) node {\rotatebox{60}{$\ddots$}};
	\filldraw (1.5,0.6) circle (0.06) node[above]{$(1)$};
	\filldraw (1.5,0.6) --++(0.5,0.2) circle (0.06) node[above]{$(1)$};
	\filldraw (2,0.8) --++(0.7,0) circle (0.06) node[above]{$(1)$};
	\draw (3.4,0.8) node{$\cdots$}; 
	\filldraw (4,0.8) circle (0.06) node[above]{$(1)$};	
	\filldraw (4,0.8) --++(0.7,0) circle (0.06) node[above]{$(1)$};
	\draw[->] (4.7,0.8) --++(0.5,0.5);
	
	\filldraw (4.7,0.8) --++(0.5,-0.4) circle (0.06) node[above]{$(1)$};
	\filldraw (5.2,0.4) --++(0.5,-0.4) circle (0.06) node[right]{$(2k-1)$};
	
	\filldraw (0,0) --++(0.5,-0.2) circle (0.06) node[below]{$(1)$};
	\draw (1,-0.3) node {\rotatebox{10}{$\ddots$}};
	\filldraw (1.5,-0.6) circle (0.06) node[below]{$(1)$};
	\filldraw (1.5,-0.6) --++(0.5,-0.2) circle (0.06) node[below]{$(1)$};
	\filldraw (2,-0.8) --++(0.7,0) circle (0.06) node[below]{$(1)$};
	\draw (3.4,-0.8) node{$\cdots$}; 
	\filldraw (4,-0.8) circle (0.06) node[below]{$(1)$};	
	\filldraw (4,-0.8) --++(0.7,0) circle (0.06) node[below]{$(1)$};
	\draw[->] (4.7,-0.8) --++(0.5,-0.5);	
	
	\filldraw (4.7,-0.8) --++(0.5,0.4) circle (0.06) node[below]{$(1)$};
	\filldraw (5.2,-0.4) --++(0.5,0.4) ;

	\draw (1.8,-0) node[right]{$\vdots$ $(2k-1)$ copies};
	
	\draw (9,0) node{($n=2k$)};
	
	\coordinate (P) at (0,-3);
		
	\draw[->] (P)++ (0,0) --++(-1,0);
	\filldraw (P)++(0,0) circle (0.06) node[above]{$(1)$};
	\filldraw (P)++(0,0) --++(0.5,0.2) circle (0.06) node[above]{$(1)$};
	\draw (P)++(1,0.5) node {\rotatebox{60}{$\ddots$}};
	\filldraw  (P)++(1.5,0.6) circle (0.06) node[above]{$(1)$};
	\filldraw (P)++(1.5,0.6) --++(0.5,0.2) circle (0.06) node[above]{$(2)$};
	\filldraw (P)++(2,0.8) --++(0.7,0) circle (0.06) node[above]{$(1)$};
	\draw (P)++(3.4,0.8) node{$\cdots$}; 
	\filldraw (P)++(4,0.8) circle (0.06) node[above]{$(1)$};	
	\filldraw (P)++(4,0.8) --++(0.7,0) circle (0.06) node[above]{$(1)$};
	\draw[->] (P)++(4.7,0.8) --++(0.5,0.5);
	
	\filldraw (P)++(4.7,0.8) --++(0.5,-0.4) circle (0.06) node[above]{$(1)$};
	\filldraw (P)++ (5.2,0.4) --++(0.5,-0.4) circle (0.06) node[right]{$(2k)$};
	
	\filldraw (P)++(0,0) --++(0.5,-0.2) circle (0.06) node[below]{$(1)$};
	\draw (P)++(1,-0.3) node {\rotatebox{10}{$\ddots$}};
	\filldraw (P)++(1.5,-0.6) circle (0.06) node[below]{$(1)$};
	\filldraw (P)++(1.5,-0.6) --++(0.5,-0.2) circle (0.06) node[below]{$(2)$};
	\filldraw (P)++(2,-0.8) --++(0.7,0) circle (0.06) node[below]{$(1)$};
	\draw (P)++(3.4,-0.8) node{$\cdots$}; 
	\filldraw (P)++(4,-0.8) circle (0.06) node[below]{$(1)$};	
	\filldraw (P)++(4,-0.8) --++(0.7,0) circle (0.06) node[below]{$(1)$};
	\draw[->] (P)++(4.7,-0.8) --++(0.5,-0.5);	
	
	\filldraw (P)++(4.7,-0.8) --++(0.5,0.4) circle (0.06) node[below]{$(1)$};
	\filldraw (P)++(5.2,-0.4) --++(0.5,0.4) ;

	\draw (P)++(1.8,-0) node[right]{$\vdots$ $2k$ copies};
	
	\draw (P)++ (9,0) node{($n=2k+1$)};

	\end{tikzpicture}
	\]
where each inserted string graph with dots has $(k-1)$ vertices if $n=2k$ and $k$ vertices if $n=2k+1$, and all non-arrow edges have signature $-$.
By applying formula (\ref{eq:mult}) and removing arrows, we have
	
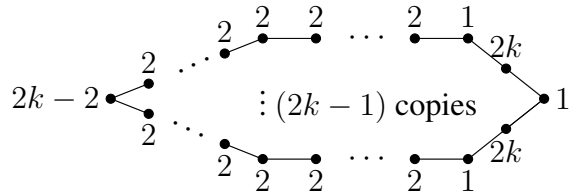
\[
	\begin{tikzpicture}
	
	\filldraw (0,0) circle (0.06) node[left]{$2k-2$};
	\filldraw (0,0) --++(0.5,0.2) circle (0.06) node[above]{$2$};
	\draw (1,0.5) node {\rotatebox{60}{$\ddots$}};
	\filldraw (1.5,0.6) circle (0.06) node[above]{$2$};
	\filldraw (1.5,0.6) --++(0.5,0.2) circle (0.06) node[above]{$2$};
	\filldraw (2,0.8) --++(0.7,0) circle (0.06) node[above]{$2$};
	\draw (3.4,0.8) node{$\cdots$}; 
	\filldraw (4,0.8) circle (0.06) node[above]{$2$};	
	\filldraw (4,0.8) --++(0.7,0) circle (0.06) node[above]{$1$};

	\filldraw (4.7,0.8) --++(0.5,-0.4) circle (0.06) node[above]{$2k$};
	\filldraw (5.2,0.4) --++(0.5,-0.4) circle (0.06) node[right]{$1$};
	
	\filldraw (0,0) --++(0.5,-0.2) circle (0.06) node[below]{$2$};
	\draw (1,-0.3) node {\rotatebox{10}{$\ddots$}};
	\filldraw (1.5,-0.6) circle (0.06) node[below]{$2$};
	\filldraw (1.5,-0.6) --++(0.5,-0.2) circle (0.06) node[below]{$2$};
	\filldraw (2,-0.8) --++(0.7,0) circle (0.06) node[below]{$2$};
	\draw (3.4,-0.8) node{$\cdots$}; 
	\filldraw (4,-0.8) circle (0.06) node[below]{$2$};	
	\filldraw (4,-0.8) --++(0.7,0) circle (0.06) node[below]{$1$};

	\filldraw (4.7,-0.8) --++(0.5,0.4) circle (0.06) node[below]{$2k$};
	\filldraw (5.2,-0.4) --++(0.5,0.4) ;

	\draw (1.8,-0) node[right]{$\vdots$ $(2k-1)$ copies};
	
	\draw (8,0) node{($n=2k$)}; 
	
	\coordinate (P) at (0,-3);
	
	\filldraw (P)++ (0,0) circle (0.06) node[left]{$2k-1$};
	\filldraw (P)++ (0,0) --++(0.5,0.2) circle (0.06) node[above]{$2$};
	\draw (P)++ (1,0.5) node {\rotatebox{60}{$\ddots$}};
	\filldraw (P)++ (1.5,0.6) circle (0.06) node[above]{$3$};
	\filldraw (P)++  (1.5,0.6) --++(0.5,0.2) circle (0.06) node[above]{$1$};
	\filldraw (P)++  (2,0.8) --++(0.7,0) circle (0.06) node[above]{$3$};
	\draw (P)++ (3.4,0.8) node{$\cdots$}; 
	\filldraw (P)++ (4,0.8) circle (0.06) node[above]{$2$};	
	\filldraw (P)++ (4,0.8) --++(0.7,0) circle (0.06) node[above]{$1$};

	\filldraw (P)++ (4.7,0.8) --++(0.5,-0.4) circle (0.06) ++(0.3,0)node[above]{$2k+1$};
	\filldraw (P)++ (5.2,0.4) --++(0.5,-0.4) circle (0.06) node[right]{$1$};
	
	\filldraw (P)++ (0,0) --++(0.5,-0.2) circle (0.06) node[below]{$2$};
	\draw(P)++  (1,-0.3) node {\rotatebox{10}{$\ddots$}};
	\filldraw (P)++ (1.5,-0.6) circle (0.06) node[below]{$3$};
	\filldraw (P)++ (1.5,-0.6) --++(0.5,-0.2) circle (0.06) node[below]{$1$};
	\filldraw (P)++ (2,-0.8) --++(0.7,0) circle (0.06) node[below]{$3$};
	\draw (P)++ (3.4,-0.8) node{$\cdots$}; 
	\filldraw(P)++  (4,-0.8) circle (0.06) node[below]{$2$};	
	\filldraw (P)++ (4,-0.8) --++(0.7,0) circle (0.06) node[below]{$1$};

	\filldraw(P)++  (4.7,-0.8) --++(0.5,0.4) circle (0.06) ++(0.3,0)node[below]{$2k+1$};
	\filldraw (P)++ (5.2,-0.4) --++(0.5,0.4) ;

	\draw (P)++ (1.8,-0) node[right]{$\vdots$ $2k$ copies};
	
	\draw (P)++ (8,0) node{($n=2k+1$)}; 
	\end{tikzpicture}
	\]
We focus on each string:
	
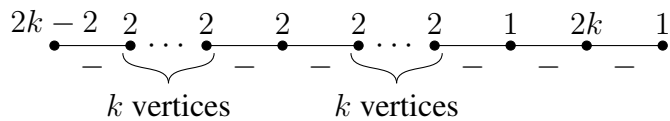
\[
	\begin{tikzpicture}
	\filldraw (0,0) circle (0.06) node[above]{$2k-2$};
	\filldraw (0,0) --++(1,0) circle (0.06) node[above]{$2$};
	\draw (1.5,0) node{$\cdots$};
	\filldraw (2,0) circle (0.06) node[above]{$2$};
	\filldraw (2,0) --++(1,0) circle (0.06) node[above]{$2$};
	\filldraw (3,0) --++(1,0) circle (0.06) node[above]{$2$};
	\draw (4.5,0) node{$\cdots$};
	\filldraw (5,0) circle (0.06) node[above]{$2$};
	\filldraw (5,0) --++(1,0) circle (0.06) node[above]{$1$};
	\filldraw (6,0) --++(1,0) circle (0.06) node[above]{$2k$};
	\filldraw (7,0) --++(1,0) circle (0.06) node[above]{$1$};
	
	\draw (0.5,0) node[below]{$-$};
	\draw (2.5,0) node[below]{$-$};
	\draw (3.5,0) node[below]{$-$};	
	\draw (5.5,0) node[below]{$-$};
	\draw (6.5,0) node[below]{$-$};
	\draw (7.5,0) node[below]{$-$};	
	
	\draw (0.9,-0.1) to[out=-90,in=90]++(0.6,-0.4) node[below]{$(k-1)$ vertices};
	\draw (2.1,-0.1) to[out=-90,in=90]++(-0.6,-0.4);
	\draw (3.9,-0.1) to[out=-90,in=90]++(0.6,-0.4) node[below]{$(k-1)$ vertices};
	\draw (5.1,-0.1) to[out=-90,in=90]++(-0.6,-0.4);
	
	\draw (10,0) node[right]{$(n=2k)$};
	 
	\coordinate (P) at (-1,-2); 
	 
	\filldraw (P) ++ (0,0) circle (0.06) node [above]{$2k-1$};
	\filldraw (P) ++ (0,0)--++(1,0) circle (0.06) node [above]{$2$};
	\draw (P) ++ (1.5,0) node {$\cdots$};
	\filldraw (P) ++ (2,0) circle (0.06) node [above]{$2$};
	\filldraw (P) ++ (2,0)--++(1,0) circle (0.06) node [above]{$3$};
	\filldraw (P) ++ (3,0)--++(1,0) circle (0.06) node [above]{$1$};
	\filldraw (P) ++ (4,0)--++(1,0) circle (0.06) node [above]{$3$};
	\filldraw (P) ++ (5,0)--++(1,0) circle (0.06) node [above]{$2$};
	\draw (P) ++ (6.5,0) node {$\cdots$};
	\filldraw (P) ++ (7,0) circle (0.06) node [above]{$2$};		
	\filldraw (P) ++ (7,0)--++(1,0) circle (0.06) node [above]{$1$};
	\filldraw (P) ++ (8,0)--++(1,0) circle (0.06) node [above]{$2k+1$};
	\filldraw (P) ++ (9,0)--++(1,0) circle (0.06) node [above]{$1$};
	
	\draw (P) ++ (0.5,0) node[below] {$-$};
	\draw (P) ++ (2.5,0) node[below] {$-$};
	\draw (P) ++ (3.5,0) node[below] {$-$};	
	\draw (P) ++ (4.5,0) node[below] {$-$};	
	\draw (P) ++ (5.5,0) node[below] {$-$};	
	\draw (P) ++ (7.5,0) node[below] {$-$};	
	\draw (P) ++ (8.5,0) node[below] {$-$};	
	\draw (P) ++ (9.5,0) node[below] {$-$};	
			
	\draw (P) ++ (0.9,-0.1) to[out=-90,in=90]++(0.6,-0.4) node[below]{$(k-1)$ vertices};
	\draw (P) ++ (2.1,-0.1) to[out=-90,in=90]++(-0.6,-0.4);
	\draw (P) ++ (5.9,-0.1) to[out=-90,in=90]++(0.6,-0.4) node[below]{$(k-1)$ vertices};
	\draw (P) ++  (7.1,-0.1) to[out=-90,in=90]++(-0.6,-0.4);
	
	\draw (P)++(11,0) node[right]{$(n=2k+1)$};
	\end{tikzpicture}
	\]
By plumbing calculus ($\pm 2$-alteration, blowing down, and $0$-chain absorption), for both case where $n$ is even or odd, each string changes to
	
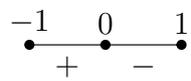
\[
	\begin{tikzpicture}
	\filldraw (0,0) circle (0.06) node[above]{$-1$};
	\filldraw (0,0) --++(1,0) circle (0.06) node[above]{$0$};
	\filldraw (1,0) --++(1,0) circle (0.06) node[above]{$1$};
	\draw (0.5,0) node[below]{$+$};
	\draw (1.5,0) node[below]{$-$};
	\end{tikzpicture}
	\]
Therefore, the plumbing graph becomes as follows:
	
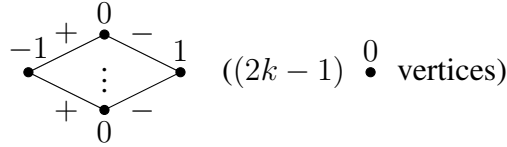
\[
	\begin{tikzpicture}
	\filldraw (0,0) circle (0.06) node[above]{$-1$};
	\filldraw (0,0) --++(1,0.5) circle (0.06) node[above]{$0$};
	\filldraw (1,0.5) --++(1,-0.5) circle (0.06) node[above]{$1$};
	\draw (1,-0.5) --++(1,0.5);
	\draw (0.5,0.8) node[below]{$+$};
	\draw (1.5,0.8) node[below]{$-$};
	\draw (1,0) node{$\vdots$};
	
	\filldraw (0,0) --++(1,-0.5) circle (0.06) node[below]{$0$};
	\draw (0.5,-0.8) node[above]{$+$};
	\draw (1.5,-0.8) node[above]{$-$};
	
	\draw (2.4,0) node[right]{($(n-1)$};
	\filldraw (4.5,0) circle (0.06) node[above]{$0$};
	\draw (4.7,0) node[right] {vertices)};
	
	\end{tikzpicture}
	\]	
Finally, by $0$-chain absorption and oriented handle absorption, we have 
	\[
	\begin{tikzpicture}
	\filldraw (0,0) circle (0.06) node [above]{$0$};
	\draw (0,0) node [below] {$\left [ n-2 \right ]$};
	\end{tikzpicture}
	\]	

Therefore, for the both case where $n$ is odd and even, $\partial \overline{F}$ is diffeomorphc to $S^1 \times \Sigma_{n-2} $, where $\Sigma_{n-2}$ is an orientable closed surface of genus $n-2$.
\end{example}

\begin{example}\label{exam:main}
Suppose that $\cA$ is a generic arrangement with $n$ hyperplanes.
The graph $\Gamma_{\cC}$ is the complete graph of $n$ vertices whose each edge is modified as follows:
	
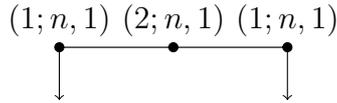
\[
	\begin{tikzpicture}
	\filldraw (0,0) circle (0.06) node[above]{$(1;n,1)$};
	\filldraw (0,0) --++(1.5,0) circle (0.06) node[above]{$(2;n,1)$};
	\filldraw (1.5,0) --++(1.5,0) circle (0.06) node[above]{$(1;n,1)$};
	\draw[->] (0,0) --++(0,-0.7);
	\draw[->] (3,0) --++(0,-0.7);
	\end{tikzpicture}
	\]
By applying the main algorithm, this string changes to
	
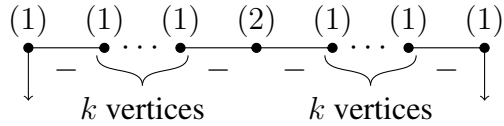
\[
	\begin{tikzpicture}
	\coordinate (P) at (0,-2);
	\filldraw (0,0) circle (0.06) node [above]{$(1)$};
	\filldraw (0,0)--++(1,0) circle (0.06) node [above]{$(1)$};
	\draw (1.5,0) node {$\cdots$};
	\filldraw (2,0) circle (0.06) node [above]{$(1)$};
	\filldraw (2,0)--++(1,0) circle (0.06) node [above]{$(1)$};
	\filldraw (3,0)--++(1,0) circle (0.06) node [above]{$(1)$};
	\draw (4.5,0) node {$\cdots$};
	\filldraw (5,0) circle (0.06) node [above]{$(1)$};		
	\filldraw (5,0)--++(1,0) circle (0.06) node [above]{$(1)$};
	
	\draw (0.5,0) node[below] {$-$};
	\draw (2.5,0) node[below] {$-$};
	\draw (3.5,0) node[below] {$-$};	
	\draw (5.5,0) node[below] {$-$};	
	
	\draw[->] (0,0) --++(0,-0.7);
	\draw[->] (6,0) --++(0,-0.7);
	\draw (0.9,-0.1) to[out=-90,in=90]++(0.6,-0.4) node[below]{$(k-1)$ vertices};
	\draw (2.1,-0.1) to[out=-90,in=90]++(-0.6,-0.4);
	\draw (3.9,-0.1) to[out=-90,in=90]++(0.6,-0.4) node[below]{$(k-1)$ vertices};
	\draw (5.1,-0.1) to[out=-90,in=90]++(-0.6,-0.4);
	
	\draw (7,0) node[right]{($n=2k$)};

	\filldraw (P)++ (0,0) circle (0.06) node [above]{$(1)$};
	\filldraw (P)++(0,0)--++(1,0) circle (0.06) node [above]{$(1)$};
	\draw (P)++(1.5,0) node {$\cdots$};
	\filldraw (P)++(2,0) circle (0.06) node [above]{$(1)$};
	\filldraw (P)++(2,0)--++(1,0) circle (0.06) node [above]{$(2)$};
	\filldraw (P)++(3,0)--++(1,0) circle (0.06) node [above]{$(1)$};
	\draw (P)++(4.5,0) node {$\cdots$};
	\filldraw (P)++(5,0) circle (0.06) node [above]{$(1)$};		
	\filldraw (P)++(5,0)--++(1,0) circle (0.06) node [above]{$(1)$};
	
	\draw (P)++(0.5,0) node[below] {$-$};
	\draw (P)++(2.5,0) node[below] {$-$};
	\draw (P)++(3.5,0) node[below] {$-$};	
	\draw (P)++(5.5,0) node[below] {$-$};	
	
	\draw[->] (P)++(0,0) --++(0,-0.7);
	\draw[->] (P)++(6,0) --++(0,-0.7);
	\draw (P)++(0.9,-0.1) to[out=-90,in=90]++(0.6,-0.4) node[below]{$k$ vertices};
	\draw (P)++(2.1,-0.1) to[out=-90,in=90]++(-0.6,-0.4);
	\draw (P)++(3.9,-0.1) to[out=-90,in=90]++(0.6,-0.4) node[below]{$k$ vertices};
	\draw (P)++(5.1,-0.1) to[out=-90,in=90]++(-0.6,-0.4);	
	
	\draw (P)++(7,0) node[right]{($n=2k+1$)};
	\end{tikzpicture}
	\]
by inserting the string $Str^{-}(1,2;2k \mid 0,0;1)$ $(n=2k)$ or $Str^{-}(1,2;2k+1 \mid 0,0;1)$ $(n=2k+1)$ with trivial covering data. For example, the complete output of the main algorithm for $n=3$ is as follows:
	
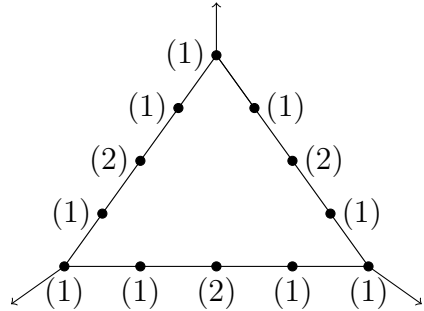
\[
	\begin{tikzpicture}
	\filldraw (0,0) circle (0.06) node [below]{$(1)$};
	\filldraw (0,0)--++(1,0) circle (0.06) node [below]{$(1)$};
	\filldraw (1,0)--++(1,0) circle (0.06) node [below]{$(2)$};
	\filldraw (2,0)--++(1,0) circle (0.06) node [below]{$(1)$};
	\filldraw (3,0)--++(1,0) circle (0.06) node [below]{$(1)$};
	
	\filldraw (0,0)--++(0.5,0.7) circle (0.06) node [left]{$(1)$};	
	\filldraw (0.5,0.7)--++(0.5,0.7) circle (0.06) node [left]{$(2)$};
	\filldraw (1,1.4)--++(0.5,0.7) circle (0.06) node [left]{$(1)$};
	\filldraw (1.5,2.1)--++(0.5,0.7) circle (0.06) node [left]{$(1)$};
	
	\filldraw (4,0)--++(-0.5,0.7) circle (0.06) node [right]{$(1)$};	
	\filldraw (3.5,0.7)--++(-0.5,0.7) circle (0.06) node [right]{$(2)$};
	\filldraw (3,1.4)--++(-0.5,0.7) circle (0.06) node [right]{$(1)$};
	\filldraw (2.5,2.1)--++(-0.5,0.7) ;	
	
	\draw[->] (0,0)--++(-0.7,-0.5);
	\draw[->] (4,0)--++(0.7,-0.5);
	\draw[->] (2,2.8)--++(0,0.7);
	\end{tikzpicture}
	\]
where all non-arrow edges have signature $-$. 
By applying formula (\ref{eq:mult}) and removing arrows, each string becomes 
	
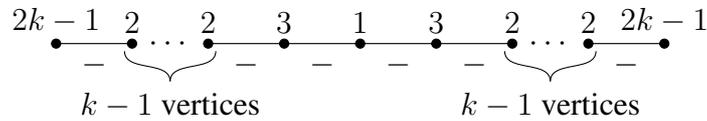
\[
	\begin{tikzpicture}
	\coordinate (P) at (-1,-2);
	
	\filldraw (0,0) circle (0.06) node [above]{$2k-2$};
	\filldraw (0,0)--++(1,0) circle (0.06) node [above]{$2$};
	\draw (1.5,0) node {$\cdots$};
	\filldraw (2,0) circle (0.06) node [above]{$2$};
	\filldraw (2,0)--++(1,0) circle (0.06) node [above]{$2$};
	\filldraw (3,0)--++(1,0) circle (0.06) node [above]{$2$};
	\draw (4.5,0) node {$\cdots$};
	\filldraw (5,0) circle (0.06) node [above]{$2$};		
	\filldraw (5,0)--++(1,0) circle (0.06) node [above]{$2k-2$};
	
	\draw (0.5,0) node[below] {$-$};
	\draw (2.5,0) node[below] {$-$};
	\draw (3.5,0) node[below] {$-$};	
	\draw (5.5,0) node[below] {$-$};	
	
	\draw (0.9,-0.1) to[out=-90,in=90]++(0.6,-0.4) node[below]{$(k-1)$ vertices};
	\draw (2.1,-0.1) to[out=-90,in=90]++(-0.6,-0.4);
	\draw (3.9,-0.1) to[out=-90,in=90]++(0.6,-0.4) node[below]{$(k-1)$ vertices};
	\draw (5.1,-0.1) to[out=-90,in=90]++(-0.6,-0.4);	

	\draw (8,0) node[right]{$(n=2k)$};

	\filldraw (P)++(0,0) circle (0.06) node [above]{$2k-1$};
	\filldraw (P)++(0,0)--++(1,0) circle (0.06) node [above]{$2$};
	\draw (P)++(1.5,0) node {$\cdots$};
	\filldraw (P)++(2,0) circle (0.06) node [above]{$2$};
	\filldraw (P)++(2,0)--++(1,0) circle (0.06) node [above]{$3$};
	\filldraw (P)++(3,0)--++(1,0) circle (0.06) node [above]{$1$};
	\filldraw (P)++(4,0)--++(1,0) circle (0.06) node [above]{$3$};
	\filldraw (P)++(5,0)--++(1,0) circle (0.06) node [above]{$2$};
	\draw (P)++(6.5,0) node {$\cdots$};
	\filldraw (P)++(7,0) circle (0.06) node [above]{$2$};		
	\filldraw (P)++(7,0)--++(1,0) circle (0.06) node [above]{$2k-1$};
	
	\draw (P)++(0.5,0) node[below] {$-$};
	\draw (P)++(2.5,0) node[below] {$-$};
	\draw (P)++(3.5,0) node[below] {$-$};	
	\draw (P)++(4.5,0) node[below] {$-$};	
	\draw (P)++(5.5,0) node[below] {$-$};	
	\draw (P)++(7.5,0) node[below] {$-$};	
	
	\draw (P)++(0.9,-0.1) to[out=-90,in=90]++(0.6,-0.4) node[below]{$(k-1)$ vertices};
	\draw (P)++(2.1,-0.1) to[out=-90,in=90]++(-0.6,-0.4);
	\draw (P)++(5.9,-0.1) to[out=-90,in=90]++(0.6,-0.4) node[below]{$(k-1)$ vertices};
	\draw (P)++(7.1,-0.1) to[out=-90,in=90]++(-0.6,-0.4);
	
	\draw (P) ++ (9,0) node[right]{$(n=2k+1)$};	
	\end{tikzpicture}
	\]
Finally, by applying plumbing calculus (blowing down, $\pm 2$-alteration, and $0$-chain absorption), for both case $n$ is even and odd, we have
	\[
	\begin{tikzpicture}
	\filldraw (0,0) circle (0.06) node [above]{$-1$};
	\filldraw (0,0)--++(1.5,0) circle (0.06) node [above]{$-n$};
	\filldraw (1.5,0)--++(1.5,0) circle (0.06) node [above]{$-1$};
	\draw (0.75,0) node[below] {$+$};
	\draw (2.25,0) node[below] {$-$};	
	\end{tikzpicture}
	\]
Thus, the resulting graph is obtained from the complete graph by adding one vertex to each edge. 
The original vertex has the Euler number $-1$, and the inserted vertex has $-n$.
Each edge of the complete graph is divided into two.
One is decorated by $-$ and the other is decorated by $+$.
\end{example}

\begin{example}\label{exam:a3}
Suppose that $\cA$ is defined by a polynomial $xyz(x-y)(y-z)(z-x)$ ($A_3$ arrangement, see Figure \ref{fig:a3}).

\begin{figure}[htbp]
\centering
\begin{tikzpicture}
\draw (-2.4,0) --(2.4,0) node[right]{$x=0$};
\draw (0,3.9) -- (0,-0.4) node[below]{$y-z=0$};
\draw (-2.2,-0.346)--(-2,0) -- (0,3.46)--++(0.2,0.346)node[right]{$y=0$};
\draw (2.2,-0.346)--(2,0) -- (0,3.46)--++(-0.2,0.346)node[left]{$z=0$};
\draw (-2.3,-0.173)--(-2,0) -- (1,1.73) --++ (0.3,0.173)node[right]{$x-y=0$};
\draw (2.3,-0.173)--(2,0) -- (-1,1.73)--++ (-0.3,0.173)node[left]{$z-x=0$};
\end{tikzpicture}
\caption{(Projectivized) $A_3$ arrangement.}
\label{fig:a3}
\end{figure}

The arrangement $\cA$ contains three double points and four triple points as the intersection points. Note that $\cA$ does not satisfy the assumption of Conjecture \ref{conj:tor} (1). 
The graph $\Gamma_{\cC}$ is described as follows:
	\[
	\begin{tikzpicture}
	\filldraw (-1,0.7) circle (0.06) node[left]{$(3;6,1)$};
	\filldraw (-1,1.4) circle (0.06) node[left]{$(3;6,1)$};
	\filldraw (-1,2.1) circle (0.06) node[left]{$(3;6,1)$};
	\filldraw (-1,2.8) circle (0.06) node[left]{$(3;6,1)$};
	\draw (-1,0.7) -- (1,0);	
	\draw (-1,0.7) -- (1,0.7);	
	\draw (-1,0.7) -- (1,2.1);	
	\draw (-1,1.4) -- (1,0);
	\draw (-1,1.4) -- (1,1.4);
	\draw (-1,1.4) -- (1,2.8);
	\draw (-1,2.1) -- (1,0.7);
	\draw (-1,2.1) -- (1,1.4);
	\draw (-1,2.1) -- (1,3.5);
	\draw (-1,2.8) -- (1,2.1);
	\draw (-1,2.8) -- (1,2.8);
	\draw (-1,2.8) -- (1,3.5);
	\draw[->] (1,3.5) --++(0.5,0);
	\draw[->] (1,2.8) --++(0.5,0);
	\draw[->] (1,2.1) --++(0.5,0);
	\draw[->] (1,1.4) --++(0.5,0);
	\draw[->] (1,0.7) --++(0.5,0);	
	\draw[->] (1,0) --++(0.5,0);
	\filldraw (1,3.5) circle (0.06) node[above]{$(1;6,1)$};
	\filldraw (1,2.8) circle (0.06);
	\filldraw (1,2.1) circle (0.06);
	\filldraw (1,1.4) circle (0.06);
	\filldraw (1,0.7) circle (0.06);
	\filldraw (1,0) circle (0.06) node[below]{$(1;6,1)$};
	\filldraw (2,1.75) circle (0.06);
	\filldraw (3,1.75) circle (0.06);
	\filldraw (4,1.75) circle (0.06) node[right]{$(2;6,1)$};
	\draw (2,1.75) -- (1,2.1);	
	\draw (2,1.75) -- (1,1.4);
	\draw (3,1.75) -- (1,2.8);
	\draw (3,1.75) -- (1,0.7);
	\draw (4,1.75) -- (1,3.5);
	\draw (4,1.75) -- (1,0);
	\end{tikzpicture}
	\]
The root of each arrow is weighted by $(1;6,1)$ and each double point vertex is weighted by $(2;6,1)$.
The following piece corresponding to the triple point 
	\[
	\begin{tikzpicture}
	\filldraw (-1,0) circle (0.06) node[left]{$(3;6,1)$};
	\filldraw (1,1) circle (0.06) node[above]{$(1;6,1)$};
	\filldraw (1,0) circle (0.06) node[above]{$(1;6,1)$};	
	\filldraw (1,-1) circle (0.06) node[above]{$(1;6,1)$};
	\draw[->] (-1,0)--(1,1)--++(1,0);
	\draw[->] (-1,0)--(1,0)--++(1,0);
	\draw[->] (-1,0)--(1,-1)--++(1,0);	
	\end{tikzpicture}
	\]
changes to the following piece by applying the main algorithm:
	\[
	\begin{tikzpicture}
	\filldraw (-1,0) circle (0.06) node[left]{$[1]$};
	\draw (-1,0) node[above]{$(1)$};
	\filldraw (1,1.5) circle (0.06) node[above]{$(1)$};
	\filldraw (1,0) circle (0.06) node[above]{$(1)$};	
	\filldraw (1,-1.5) circle (0.06) node[above]{$(1)$};
	\filldraw (-0,0.75) circle (0.06) node[above]{$(1)$}; 
	\filldraw (-0,0) circle (0.06) node[above]{$(1)$};
	\filldraw (-0,-0.75) circle (0.06) node[above]{$(1)$};		
	\draw[->] (-1,0)--(1,1.5)--++(1,0);
	\draw[->] (-1,0)--(1,0)--++(1,0);
	\draw[->] (-1,0)--(1,-1.5)--++(1,0);
	\draw (-0.4,0.35) node{$-$};
	\draw (-0.4,-0.1) node{$-$};
	\draw (-0.4,-0.65) node{$-$};
	\draw (0.5,0.95) node{$-$};
	\draw (0.5,-0.1) node{$-$};
	\draw (0.5,-1.3) node{$-$};		
	\end{tikzpicture}
	\]
Here, we inserted the string 
$Str^{-} (1;3,6 \mid 0,0;1)
= \begin{tikzpicture}
\draw[<->] (-1,0) -- (1,0);
\filldraw (0,0) circle(0.06) node[above]{$(1)$};
\draw (1,0) node[above]{$(1)$};
\draw (-1,0) node[above]{$(1)$};
\end{tikzpicture}$ 
with trivial covering data.
By applying formula (1) and removing arrows, this piece becomes
	\[
	\begin{tikzpicture}
	\filldraw (-1,0) circle (0.06) node[left]{$[1]$};
	\draw (-1,0) node[above]{$3$};
	\filldraw (1,1.5) circle (0.06) node[above]{$-2$};
	\filldraw (1,0) circle (0.06) node[above]{$-2$};	
	\filldraw (1,-1.5) circle (0.06) node[above]{$-2$};
	\filldraw (-0,0.75) circle (0.06) node[above]{$2$}; 
	\filldraw (-0,0) circle (0.06) node[above]{$2$};
	\filldraw (-0,-0.75) circle (0.06) node[above]{$2$};		
	\draw (-1,0)--(1,1.5);
	\draw (-1,0)--(1,0);
	\draw (-1,0)--(1,-1.5);
	\draw (-0.4,0.35) node{$-$};
	\draw (-0.4,-0.1) node{$-$};
	\draw (-0.4,-0.65) node{$-$};
	\draw (0.5,0.95) node{$-$};
	\draw (0.5,-0.1) node{$-$};
	\draw (0.5,-1.3) node{$-$};	
	\end{tikzpicture}
	\]
By performing plumbing calculus ($\pm 2$-alteration), this piece changes to
\[
	\begin{tikzpicture}
	\filldraw (-1,0) circle (0.06) node[left]{$[1]$};
	\draw (-1,0) node[above]{$0$};
	\filldraw (1,1.5) circle (0.06) node[above]{$-1$};
	\filldraw (1,0) circle (0.06) node[above]{$-1$};	
	\filldraw (1,-1.5) circle (0.06) node[above]{$-1$};
	\filldraw (-0,0.75) circle (0.06) node[above]{$-2$}; 
	\filldraw (-0,0) circle (0.06) node[above]{$-2$};
	\filldraw (-0,-0.75) circle (0.06) node[above]{$-2$};		
	\draw (-1,0)--(1,1.5);
	\draw (-1,0)--(1,0);
	\draw (-1,0)--(1,-1.5);
	\draw (-0.4,0.35) node{$-$};
	\draw (-0.4,-0.1) node{$-$};
	\draw (-0.4,-0.65) node{$-$};
	\draw (0.5,0.95) node{$+$};
	\draw (0.5,-0.1) node{$+$};
	\draw (0.5,-1.3) node{$+$};	
	\end{tikzpicture}
	\]
On the other hand, similarly to generic arrangements, each string corresponding to double points changes to the following string:
	\[
	\begin{tikzpicture}
	\filldraw (0,0) circle (0.06) node [above]{$-1$};
	\filldraw (0,0)--++(1.5,0) circle (0.06) node [above]{$-6$};
	\filldraw (1.5,0)--++(1.5,0) circle (0.06) node [above]{$-1$};
	\draw (0.75,0) node[below] {$+$};
	\draw (2.25,0) node[below] {$-$};	
	\end{tikzpicture}
	\]
Note that the above two plumbing calculus are performed simultaneously.
Combining these pieces, we obtain a plumbing graph for the Milnor fiber boundary of the $A_3$ arrangement. 
Using Theorem \ref{thm:homology}, a direct computation show us that 
\[
H_{1} (\partial F; \bZ) \cong \bZ^{19} \oplus \bZ^{2}_{2}.
\]
Therefore, this gives an example such that $H_{1} (\partial F; \bZ)$ consists a direct summand other than $\bZ_{n}$, where $n=|\cA|$.
It would be interesting to investigate the relationship between the order of the torsion in $H_{1} (\partial F; \bZ)$ and the number $n/\gcd(m_{j}, n)$ for each intersection point.
\end{example}

\end{document}